\documentclass[12pt]{article}
\usepackage{amsfonts}
\usepackage{stmaryrd}
\usepackage{amssymb}
\usepackage{amsthm}
\usepackage{latexsym,amsmath,amscd}
\usepackage{color}
\usepackage{dsfont}

\usepackage{amssymb}
\usepackage{latexsym}
\usepackage{cleveref}

\newtheorem{theorem}{Theorem}[section]
\newtheorem{lem}[theorem]{Lemma}
\newtheorem{pro}[theorem]{Proposition}
\newtheorem{cor}[theorem]{Corollary}
\newtheorem{rem}[theorem]{Remark}
\newtheorem{rems}[theorem]{Remarks}
\newtheorem{ex}[theorem]{Example}
\newtheorem{defi}[theorem]{Definition}
\newtheorem{hyp}[theorem]{Assumption}

\crefname{theorem}{Theorem}{theorem}
\crefname{lem}{Lemma}{lemmas}
\crefname{pro}{Proposition}{proposition}
\crefname{cor}{Corollary}{corollary}
\crefname{rem}{Remark}{remark}
\crefname{rems}{Remarks}{remarks}
\crefname{ex}{Remarks}{remarks}
\crefname{defi}{Definition}{definition}
\crefname{hyp}{Assumption}{assumptions}
\crefname{con}{Conjecture}{Conjecture}
\crefname{equation}{}{Eqs.}

\newcommand{\be}{\begin{equation}}
\newcommand{\ee}{\end{equation}}
\newcommand{\bde}{\begin{displaymath}}
\newcommand{\ede}{\end{displaymath}}
\newcommand{\beq}{\begin{eqnarray*}}
\newcommand{\eeq}{\end{eqnarray*}}
\newcommand{\beqa}{\begin{eqnarray}}
\newcommand{\eeqa}{\end{eqnarray}}
\newcommand{\bel }{\left\{\begin{array}{ll}}
\newcommand{\eel}{\cr \end{array} \right.}

\newcommand{\cE}[3]{\mathbb{E}_{#1}\left(\left. #2\, \right| #3 \right)}
\newcommand{\cP}[3]{#1\left(\left. #2\, \right| #3 \right)}
\newcommand{\seq}[1]{{\lbrace #1 \rbrace}}

\newcommand{\pbracket}[2]{\left< #1, #2 \right>}
\newcommand{\expect}[2]{\mathbb{E}_{#1}( #2 ) }

\newcommand{\dcb}{\begin{array}{lll}}
\newcommand{\dce}{\end{array}}
\newcommand{\ebe}{\begin{enumerate}\setlength{\baselineskip}{13pt}\setlength{\parskip}{0pt}}
\newcommand{\dbe}{\end{enumerate}}

\def \ind{1\!\!1\!}
\def\cro#1{\langle #1\rangle}

\def\wt{\widetilde}
\def\wh{\widehat }
\def\F{{\cal F}}

\def\G{{\cal G}}
\def\H{{\cal H}}
\def\K{{\cal K}}

\def\ff{{\mathbb F}}
\def\gg{{\mathbb G}}
\def\hh{{\mathbb H}}
\def\kk{{\mathbb K}}
\def\rr{{\mathbb R}}

\def\nn{{\mathbb N}}
\def\Q{{\mathbb Q}}
\def\P{{\mathbb P}}

\def\I{\mathds{1}}

\newcommand{\bt}{\begin{theorem}}
\newcommand{\et}{\end{theorem}}
\newcommand{\bl}{\begin{lem}}
\newcommand{\el}{\end{lem}}
\newcommand{\bp}{\begin{pro}}
\newcommand{\ep}{\end{pro}}
\newcommand{\bcor}{\begin{cor}}
\newcommand{\ecor}{\end{cor}}

\newcommand{\bd}{\begin{defi} \rm }
\newcommand{\ed}{\end{defi}}
\newcommand{\brem }{\begin{rem} \rm }
\newcommand{\erem }{\end{rem}}
\newcommand{\brems }{\begin{rems} \rm }
\newcommand{\erems }{\end{rems}}
\newcommand{\bhyp }{\begin{hyp} \rm }
\newcommand{\ehyp }{\end{hyp}}
\newcommand{\bex}{\begin{ex} \rm }
\newcommand{\eex}{\end{ex}}

\newcommand{\llb}{\llbracket\,}
\newcommand{\rrb}{\,\rrbracket}

\title{{\Large \bf An enlargement  of filtration formula with application to progressive enlargement with multiple  random times\footnote{Preliminary version}}\vskip 75 pt}

\author{Monique Jeanblanc, Shiqi Song
\thanks{This research   benefited from the support of the
`Chaire Risque de cr\'edit', F\'ed\'eration Bancaire Fran\c caise
}\\
Laboratoire Analyse et Probabilit\'es\\
Universit\'e d'\'Evry Val d'Essonne  \\
91037 \'Evry Cedex, France\\\\
Libo, Li\footnotemark[2]\,\,\,\thanks{This research benefited from the support of Japan Science and Technology Agency}\\
Department of Mathematical Sciences\\
Ritsumeikan University\\
525-8577 Shiga, Japan
}
\date{}

\newcommand{\dcbs}{\begin{array}{lll}}
\newcommand{\dces}{\end{array}}
\newcommand{\ebes}{\begin{enumerate}\setlength{\baselineskip}{13pt}\setlength{\parskip}{5pt}}
\newcommand{\dbes}{\end{enumerate}}

\newcommand{\btau}{{\boldsymbol{\tau}}}

\begin{document}

\maketitle

\vskip 140 pt

\pagebreak

\begin{abstract}

Given a reference  filtration $\ff$, we develop in this work  a
generic method for computing the semimartingale decomposition of
$\ff$-martingales in some specific enlargements of  $\ff$. This method
  is then applied to the study of progressive enlargement
with multiple non-ordered random times, for which explicit
decompositions can be obtained under the absolute continuity
condition of Jacod.
\end{abstract}


\vfill \eject

\tableofcontents
 \vfill\eject

\section{Introduction}

In this paper, we work on a filtered probability space $(\Omega,
\F, \P)$ endowed with a filtration $\ff=(\F_t)_{t\geq 0}$
satisfying the usual conditions with $\F_\infty \subset \F$. For a semi-martingale $X$
and a predictable process $H$, we denote by $H\centerdot X$ the stochastic integral of $H$ with respect to $X$, whenever it is well defined.
 The optional (resp. predictable) $\sigma$-algebra
generated by a filtration $\ff$ is denoted by $\mathcal{O}(\ff)$
(resp. $\mathcal{P}(\ff)$). For the ease of language, for any
$\ff$-special semimartingale $X$,
 the $\ff$-predictable
process of finite variation  $A$ in its $\ff$-semimartingale
decomposition $X=M+A$ is called the $\ff$-drift of $X$.

Given a reference filtration $\ff$ and a filtration $\gg$ such
that $\ff \subset \gg$, one aim in the theory of enlargement is to
study whether the hypothesis $(H')$ is satisfied between $\ff$ and
$\gg$, i.e., whether any $\ff$-martingale is a
$\gg$-semi-martingale. Traditionally, one attempts this directly
by looking only at the filtrations $\ff$ and $\gg$. In the current
literature, this direct approach has been used in the study of
progressive (initial) enlargement of $\ff$ with a random time (a
 non negative random variable), including (but not limited to) the works of
Jeulin \cite{J1,J2}, Jeulin and Yor \cite{JY1} and Jacod
\cite{JJ}.

However, it is often difficult to study the hypothesis
$(H')$ directly between $\ff$ and $\gg$, and, in some cases one
can take advantage of the following result of Stricker
\cite{St}.

\bp\label[pro]{fgk} Let $\ff$, $\gg$ and $\wh \ff$  be filtrations
such that $\ff \subset \gg \subset \wh \ff$. If the hypothesis
$(H')$ is satisfied between $\ff$ and $\wh \ff$, then the
hypothesis $(H')$ is also satisfied between $\ff$ and $\gg$. \ep

 In
addition, if the $\wh \ff$-semimartingale decomposition of
$\ff$-martingales is known then (at least theoretically) one can
make use of the $\gg$-optional and $\gg$-dual predictable
projections to find the $\gg$-semimartingale decomposition  of
$\ff$-martingales.

In the current literature, Proposition \ref{fgk} has been used in
Jeanblanc and Le Cam \cite{JL}, Callegaro et al. \cite{CJZ} and
Kchia et al. \cite{KLP}, to study the relationship between the
filtrations $\ff \subset \ff^\tau \subset \gg^\tau$, where
$\ff^\tau$ (resp. $\gg^\tau$)  is  the progressive   (resp.
initial) enlargement of $\ff$ with the random time $\tau$.

If one assumes that the
$\gg^\tau$-semimartingale decomposition of $\ff$ martingales is
known, then essentially by exploiting the property that for any
$\gg^\tau$-predictable process $V^*$, the process $(V^*-V^*_{\tau})\I_{\rrb \tau,
\infty \llb}$ is $\ff^\tau$-predictable, one can derive the $\ff^\tau$-semimartingale decomposition of
$\ff$-martingales without calculating any $\ff^\tau$-dual predictable projection (see Lemma 2 and Theorem 3 in \cite{KLP}).

In this paper, we also take advantage of the above Proposition to study
the hypothesis $(H')$ between $\ff$ and any enlargement $\gg$ of $\ff$
that satisfies a specific structure (see \cref{g}).  The filtration $\gg$ arises naturally while
studying the hypothesis $(H')$ between $\ff$ and its progressive
enlargement with minimum and maximum of two random times. However,
we must stress that, although, we apply our result to the study of progressive enlargement with random times, our setup is different from 
the usual progressive enlargement framework, and the progressive enlargement with a single random time cannot be retrieved from our setting as a specific case. Therefore, unlike the case
studied in \cite{CJZ}, \cite{JL} and \cite{KLP}, one cannot
exploit the specific structure between {\it progressive} and {\it
initial} enlargement with a single random time in computing the
$\gg$-semimartingale decomposition. (for details, see \cref{hh1}.)

The present work is divided into two parts. In \cref{sH1}, a
generic technique for computing the semimartingale decomposition
of $\ff$-martingales in an enlargement of $\ff$ satisfying a
specific structure (namely a filtration $\gg$ satisfying \cref{g})
is developed. We construct a filtration $\wh \ff$ ({\it the direct
sum filtration}) that satisfies $\gg\subset \wh \ff$  and we show
that the hypothesis $(H')$ holds between $\ff$ and $\gg$ by
showing that the hypothesis $(H')$ holds between $\ff$ and $\wh
\ff$. The explicit $\gg$-decomposition formula is fostered in
\cref{gde}, which is obtained from projecting the $\wh
\ff$-semimartingale decomposition computed in \cref{fhatdecomp}. In
\cref{applied}, we apply the results
of \cref{sH1} to the study of progressive
enlargement of $\ff$ with multiple non-ordered random times. Under
the additional assumption that the joint $\ff$-conditional
distribution of the given family of random times is absolutely
continuous with respect to a non-atomic measure,   the
semimartingale decomposition is fully explicit (under mild
integrability conditions). Progressive enlargement with multiple
 ordered random times have also been studied in \cite{GN} and
\cite{JKP}, where the authors assume that the random times are
ordered with random marks. Our approach is different from the method
presented in \cite{GN} and \cite{JKP}, as it is generic and does not depend on the
techniques developed in the literature on progressive enlargement with random times.

%

\section{The setup and the main results}\label{sH1}
 We work on a  probability space $(\Omega, \F,\P)$ endowed with a
filtration $\ff= (\F_t)_{t\geq 0}$ satisfying the usual
conditions.  Keeping \cref{fgk} in mind, for a filtration $\gg$
that satisfies some specific assumptions, we construct in
\cref{sH12}, a filtration $\wh \ff$ such that $\ff \subset
\gg\subset \wh\ff$ and that the hypothesis $(H')$ is satisfied
between $\ff$ and $\wh \ff$. In   \cref{fhatdecomp}, we give
the $\wh\ff$-semimartingale decomposition of $\ff$-local martingales.
Then in \cref{hh1} the $\gg$-semimartingale decomposition of
$\ff$-martingales is deduced from \cref{fhatdecomp}.

Before proceeding, let us introduce the  following notion, which
will be important through out this section.

\bd\label[defi]{coin} Given two $\sigma$-algebras $\K$, $\H$ and a
set $D$, we write $\K\cap D \subset \H\cap D$, if for every
$\K$-measurable set $X$, there exists a $\H$-measurable set $Y$
such that $X\cap D = Y\cap D$. We say that the $\sigma$-algebras
$\K$ and $\H$ coincide on $D$, if $\K\cap D \subset \H\cap D$  and
$\K\cap D \supset \H\cap D$, in which case, we write $\K\cap D =
\H\cap D$. \ed

Let us illustrate this notion in the classical cases of the initial
enlargement $\gg^{\tau} = (\G^{\tau}_t)_{t\geq 0}$ of $\ff$ with
the random variable $\tau$ (this  is the smallest filtration
containing $\ff$ which satisfies the usual conditions such that
$\tau$ is $\G^\tau_0$-measurable) and  of the   progressive
enlargement $\ff^{\tau} = (\F^{\tau}_t)_{t\geq 0}$ of $\ff$ with
the random time $\tau$ (the smallest filtration containing $\ff$
which satisfies the usual conditions such that $\tau$ is a
stopping times). In the   case of progressive enlargement with a
random time, we know from Jeulin \cite{J2} that for every $t$, the
$\sigma$-algebras $\G^\tau_t$ and $\F^\tau_t$ coincide on the set
$\seq{\tau \leq t}$ and that the $\sigma$-algebras $\F_t$ and
$\F^\tau_t$ coincide on the set $\seq{\tau > t}$. This fact
implies also that for any $\mathcal{P}(\ff^\tau)$-measurable
random variable (process) $V$ there exists a random variable
(process) $K$ (resp. $J$) which is $\mathcal{P}(\ff)$ (resp.
$\mathcal{P}(\gg^\tau)$) measurable
  such that \bde V\I_{\llb 0, \tau \rrb} = K\I_{\llb 0,
\tau \rrb} \quad \mathrm{and} \quad V\I_{\rrb \tau ,\infty\llb} =
J\I_{\rrb\tau,\infty \llb}. \ede

In other terms,  the $\sigma$-algebras
$\mathcal{P}(\ff^\tau)$ and $\mathcal{P}(\ff)$ coincide on the set
${\llb 0, \tau \rrb}$, while $\mathcal{P}(\ff^\tau)$ and
$\mathcal{P}(\gg^\tau)$ coincide on ${\rrb\tau,\infty \llb}$.


In this paper, we study the hypothesis $(H')$ between $\ff$ and an
enlargement of $\ff$ denoted by $\gg = (\G_t)_{t\geq 0}$ which
satisfies the following assumption. \bhyp\label[hyp]{g} The
filtration $\gg$ is such that there exists an $\F$-measurable
partition of $\Omega$ given by $\seq{D_1,\dots, D_k}$ and a family
of right-continuous filtrations $\seq{\ff^1,\dots, \ff^k}$ where
for every $i=1,\dots, k$\hfill\break (i)   $\ff\subset \ff^i$ and
$\F_\infty^i \subset \F$,\hfill\break (ii) for all $t\geq 0$, the
$\sigma$-algebras $\G_t$ and   $\F^i_t$ coincide on $D_i$.\ehyp

{In such a case, we shall say that
$(\ff,(\ff^i)_{i=1,\cdots,k}, \gg)$ satisfies Assumption \ref{g}
 with respect to the partition $(D_1,\cdots, D_k)$.}


The setting here is different from that of progressive enlargement studied in \cite{CJZ}, \cite{JL} and \cite{KLP}, and as mentioned above, we cannot retrieve from our
 framework, the progressive enlargement with a single random time. This is because, in the case of progressive enlargement with a single random time $\tau$, for every $t\geq 0$, the space $\Omega$ is partitioned into $\seq{\tau> t}$ and $\seq{\tau\leq t}$, which are time dependent. Whereas, 
in our setting, by using \cref{g}, we essentially partition the space $\Omega$ using $D_i$ for $i =1,\dots, k$, which are independent of time. In other words, unlike the progressive enlargement case, we partition the product space $\Omega\times [0,\infty[$ only in {\it `space'}, rather 
than in both {\it `space'} and {\it `time'}.

On the other hand, although possible, we do not try to generalize the idea of partitioning in both {\it `space'} and {\it `time'}, to include the well studied initial and progressive enlargement setting (cf. Kchia and Protter \cite{KP}). Our purpose is different and the fact that we do not partition in {\it `time'} is used later in the computation of $\gg$-semimartingale decomposition.

Given a filtration $\gg$ satisfying \cref{g}, the goal here is to study the $\gg$-semimartingale decomposition of an $\ff$-martingale $M$, which for every $i= 1,
\dots, k$ has $\ff^i$-semimartingale decomposition given by $M = M^i + K^i$,
where $M^i$ is an $\ff^i$-local martingale and $K^i$ is an
$\ff^i$-predictable process of finite variation. The main idea of the paper is to construct from $\gg$, a filtration $\wh \ff$ (see \eqref{whf}) such that we have $\ff\subset \gg\subset \wh \ff$, and compute the $\gg$-semimartingale decomposition of a $\ff$-martingale from its $\wh \ff$-semimartingale decomposition.

Before going into the technical details in the rest of the paper, we describe first our main results. Instead of working with the filtration $\gg$ directly, the $\wh \ff$-semimartingale decomposition of $M$ is easier to compute.

\bt\label{fhatdecomp} 
If for every $i=1,\dots, k$, the
$\ff$-martingale $M$ is an $\ff^i$-semimartingale with 
$\ff^i$-semimartingale decomposition given by $M = M^i + K^i$,
where $M^i$ is an $\ff^i$-local martingale and $K^i$ is an
$\ff^i$-predictable process of finite variation, then under
Assumption  \ref{g}, \bde \label{d1} M -
\sum_{i=1}^k\left(K^i\I_{D_i} +
\frac{\I_{D_i}}{N^i_-}\centerdot\pbracket{N^i}{M}^i\right) \ede is
an $\wh \ff$-local martingale. \et

Then we can retrieve the $\gg$-semimartingale decomposition of $M$ by calculating the $\gg$-optional projection of the above formula. We introduce, for every $i= 1,\dots, k$, 
\begin{align*}
\wt N^i &:= \,^{o,\gg}(\I_{D_i}), \quad N^i := \,^{o,\ff^i}(\I_{D_i})\\
\wh V^i &:= K^i+ \frac{1}{N^i_-}\centerdot\pbracket{N^i}{M}^i
 \end{align*}
where $\pbracket{\cdot}{\cdot}^i$ is the $\ff^i$-predictable bracket.

\bt Let  $M$ be an $\ff$-martingale, such that for every
$i= 1,\dots k$, $M= M^i + K^i$, where $M^i$ is an $\ff^i$-local
martingale and $K^i$ an $\ff^i$-predictable process of finite
variation. Then under \cref{g}, for every $i = 1,\dots, k$, the
process $\wt N^i$ belongs to
$\mathcal{L}^1(\psi(\wh V^i))$ and
\begin{align*}
M -  \sum^k_{i=1} \wt N^i_-\ast \psi(\wh V^i)
\end{align*}
is a $\gg$-local martingale. The operations $\psi$ and $\ast$ are defined below in \cref{vv} and \eqref{psivv} respectively.
\et
\brem
We show in \cref{gdel} that for every $i=1,\dots, k$, the process $\wt N^i_-\ast \psi(\wh V^i)$ is in fact equal to $(\wh V^i\I_D)^{p,\gg}$. Furthermore, without going into all the details, we point out that the operations $\psi$ and $\ast$ are essentially introduced for technical reasons, namely, to deal with the cases where $\wt N^i= \,^{o,\gg}(\I_{D_i})$ takes value zero. If one considers only $\wt N^i_-\ast \psi(\wh V^i)$ up to the first time $\wt N^i$ hits zero, then $\wt N^i_-\ast \psi(\wh V^i)$ is simply the stochastic integral of $\wt N^i$ against $\psi(\wh V^i) = \frac{\,^{p,\gg}(\wh V^i\I_D)}{\,^{p,\gg}(\I_D)}$.
\erem


\subsection{Direct Sum Filtration}\label{sH12}We first construct a
filtration $\wh \ff$, called   the direct sum filtration, such
that $\gg\subset \wh \ff$ and study its properties.   Let us
define for every $t\geq 0$ the following family of sets \be
\label{whf} \widehat\F_t := \seq{A\in \F\, |\,\forall i,\,\exists
A^{i}_t\in \F^i_t \,\,\mathrm{such\,that}\, \, A\cap D_i = A^i_t
\cap D_i}\,. \ee 
The family $\widehat \ff=(\widehat\F_t )_{t\geq
0}$ will be shown in \cref{minF} to be a right continuous
filtration such that $\ff \subset \gg \subset \wh \ff$ holds.   In
general, the inclusion $\gg \subset \wh \ff$ is strict as for $i =
1,\dots, k$, the sets $D_i$ are $\wh \F_0$-measurable (hence
$\widehat\F_t$-measurable) but not necessarily $\G_0$-measurable.

The constructed filtration $\wh \ff$ and the subsequent
$\wh \ff$-semimartingale decomposition formula derived in \cref{fhatdecomp} can be viewed as an extension of the study of initial enlargement in Meyer \cite{M5} and Yor
\cite{Y3}. In \cite{M5} and \cite{Y3}, the authors enlarged the base filtration $\ff$ with a finite valued random variable $X = \sum^k_{i=1} a_i \I_{D_i}$, where for all $i= 1,\dots, k$, $a_i \in \rr$ and $D_i =\seq{X=a_i}$. If we construct the filtration $\wh \ff$, by taking the partition of $\Omega$ to be $(\seq{X = a_i})_{i=1,\dots ,k}$ and setting for all $i = 1,\dots, k$, $\ff^i= \ff$, then we can recover from \cref{decompminmax} the semimartingale
decomposition result given in Meyer \cite{M5} and in Yor
\cite{Y3}. 

\bl \label[lem]{cond} For every $t\geq 0$ and $i = 1,\dots ,k$,
\hfill\break (i) the inclusion $D_i \subseteq
\seq{\cP{\P}{D_i}{\F^i_t} > 0}$ holds $\P$-a.s. \hfill\break (ii)
For any $\P$-integrable random variable  $\eta$, one has
  \be
\cE{\P}{\eta\I_{D_i}}{\wh \F_t} =
\I_{D_i}\frac{\cE{\P}{\eta\I_{D_i}}{\F^i_t}}{\cP{\P}{D_i}{\F^i_t}}.
\label{condeq} \ee
\begin{proof}
Let $t\geq 0$  be fixed and $i = 1,\dots ,k$. \hfill\break (i) For
$\Delta = \seq{\cP{\P}{D_i}{\F^i_t} > 0}$,  one has
$\expect{\P}{\I_{D_i}\I_{\Delta^c}} = 0$, which implies that
$\P$-a.s. $D_i \subset \Delta$.\hfill\break
(ii) For $B \in \wh \F_t$,   by definition,  there exists a set
$B^i \in \F^i_t$ such that $B\cap D_i=B^i \cap D_i$. Then we have
\begin{align*}
\expect{\P}{\eta\I_{D_i}\I_B} & = \expect{\P}{\eta\I_{D_i}\I_{B^i}} = \expect{\P}{\I_{D_i}\frac{\cE{\P}{\eta\I_{D_i}}{\F^i_t}}{\cP{\P}{D_i}{\F^i_t}} \I_{B^i}}  = \expect{\P}{\I_{D_i}\frac{\cE{\P}{\eta \I_{D_i}}{\F^i_t}}{\cP{\P}{D_i}{\F^i_t}} \I_{B}}.
\end{align*}
\end{proof}
\el

\bl\label[lem]{minF}
The family $\wh \ff = (\widehat\F_t)_{t\geq
0}$ is a right-continuous filtration and   $\gg \subset \wh \ff$.
\begin{proof}
It is not hard to check that the family $\wh \ff$ is a filtration and  that $\gg$ is a subfiltration of $\wh \ff$, therefore we only prove that $\wh \ff$ is right continuous.\hfill\break\newline
To do that, we show for a fixed $t\geq 0$, that for every set $B \in \cap_{s>t}\wh\F_s$, and for every $i = 1,\dots, k$, there exists $B_i \in
\F^i_t$ such that $ B \cap  {D_i} =  {B_i}\cap {D_i}$. The set $B$ is $\wh \F_q$ measurable for all rational number $q$ strictly greater than $t$, thus, for each rational $q>t$ and each $i=1,\dots,k$, there exists an $\F^i_q$-measurable set $B_{i,q}$, such that $B\cap D_i= B_{i,q}\cap D_i$. It is sufficient to set
\bde
B_i:= \bigcap_{n\geq 0}\bigcup_{q\in (t,t+\frac{1}{n}]} B_{i,q},
\ede 
which is $\cap_{s>t} \F^i_s$ measurable and therefore $\F^i_t$-measurable by right continuity of $\ff^i$.
\end{proof}
\el

\bl\label[lem]{pp}
Under Assumption  \ref{g}, for every $i = 1,\dots,k$, the
$\sigma$-algebras $\mathcal{P}(\wh \ff)$, $\mathcal{P}(\ff^i)$ and $\mathcal{P}(\gg)$
coincide on $D_i$.

\begin{proof}
The fact that the $\sigma$-algebras $\mathcal{P}(\wh \ff)$ and
$\mathcal{P}(\ff^i)$ coincide on $D_i$ is a straightforward
consequence of the definition of $\wh \ff$. On the other hand,
$\mathcal{P}(\ff^i)$ and $\mathcal{P}(\gg)$ coincide on $D_i$ due
to (ii) of \cref{g}.
\end{proof}

\el

%

We are now in the position to compute the $\wh \ff$-semimartingale
decomposition of $\ff$-martingales. Let us first introduce some
useful notation. For each $i = 1,\dots ,k$, we define an $\ff ^i$
local martingale as $N^i := \,^{o,\ff^i}(\I_{D_i})$, where the
process $^{o,\ff^i}(\I_{D_i})$ is the c\`adl\`ag version of the
$\ff^i$-optional projection of the random variable $\I_{D_i}$. The processes $N^i$ are
bounded and therefore are $BMO$-martingales in the filtration
$\ff^i$. 

\bp\label[pro]{decompminmax} Under Assumption  \ref{g},  if  $M^i$
is an $\ff^i$-local martingale, then \bde \wh M^i:=\left(M^i-
\frac{\I_{D_i}}{N^i_-}\centerdot
\pbracket{N^i}{M^i}^i\right)\I_{D_i} \ede is an $\wh \ff$-local
martingale.

\begin{proof}
We start by noticing that the $\ff^i$-predictable bracket
$\pbracket{M^i}{N^{i}}^i$ exists, since $N^i$ is a
$BMO$-martingale. Let $(T_n)_{n \in \nn^+}$ be a sequence of
$\ff^i$-stopping times which increases to infinity, such that, for
each $n$, the process $(\pbracket{M^i}{N^{i}}^i)^{T_n}$ is of
integrable variation and that $(M^i)^{T_n}$ is a uniformly
integrable $\ff^i$-martingale. For every $n\in \nn^+$, we
construct a sequence $r_n$ of $ \ff^i$-stopping times by setting
$r_n := \inf\seq{t
>0: N^i_{t}\leq \frac{1}{n}}$ and define
$S_{n,D_i} = (r_n\wedge T_{n})\I_{D_i} + \infty \I_{D^c_i}$, which
is a sequence of $\wh \ff$-stopping times such that $S_{n,D_i}
\rightarrow \infty$ as $n \rightarrow \infty$.

For every $n\in \nn^+$ and any bounded elementary $\wh \ff$-predictable process $\wh \xi$,
\begin{align*}
\expect{\P}{(\I_{D_i}\wh\xi\centerdot (M^i)^{S_{n,D_i}})_\infty} &=  \expect{\P}{(\I_{D_i}\xi^i\centerdot (M^i)^{r_{n}\wedge T_n} )_\infty}   \\
                                  &=  \expect{\P}{(\xi^i\I_{\llb 0,r_{n}\wedge T_n\rrb}\centerdot \pbracket{N^i}{M^i}^i)_\infty },
\end{align*}
where $\xi^i$ is the bounded $\ff^i$-predictable process associated with
the process $\wh \xi$ and the set $D_i$ by \cref{pp}, and the second equality holds by
integration by parts formula in $\ff^i$. Next, we see that
\begin{align*}
\expect{\P}{(\xi^i\I_{\llb 0,r_{n}\wedge T_n\rrb}\centerdot \pbracket{N^i}{M^i}^i)_\infty} & =  \expect{\P}{(\I_{D_i}\xi^i(N^i_-)^{-1}\I_{\llb 0,r_{n}\wedge T_n\rrb}\centerdot \pbracket{N^i}{M^{i}}^i)_\infty } \\
& =  \expect{\P}{\I_{D_i}(\wh\xi\I_{D_i}(N^i_-)^{-1}\I_{\llb 0,S_{n,D_i}\rrb}\centerdot
\pbracket{N^i}{M^{i}}^i)_\infty }
\end{align*}
where the first equality holds by taking the $\ff^i$-dual
predictable projection. This shows that \bde
\left(M^{i}-
\frac{\I_{D_i}}{N^i_-}\centerdot\pbracket{N^i}{M^{i}}^i\right)\I_{D_i} \ede is an $\wh \ff$-local martingale.
\end{proof}

\ep



Theorem \ref{fhatdecomp} is now an immediate consequence of Proposition \ref{decompminmax}. From that, we see that if the hypothesis $(H^\prime)$ is
satisfied between $\ff$ and $\ff^i$ for any $i$, then
the hypothesis $(H^\prime)$ is satisfied between $\ff$ and $\wh \ff$.

\subsection{Computation of $\gg$-semimartingale decomposition}\label{hh1}

Before proceeding with the computation of the
$\gg$-semimartingale decomposition of $\ff$-martingales, we first
summarize the current study and make some comparisons with the well
studied cases of initial and progressive enlargement. Then we
explain why our study is different from what has been previously
done in the literature.

In our setting, we have the following relationships \bde \ff
\subset \gg \subset \wh \ff \quad and \quad \ff \subset \ff^i,
\quad \forall i = 1,\dots, k \ede 
whereas, in the classic initial and progressive enlargement setting, we have $\ff \subset \ff^\tau \subset \gg^\tau$.

Similar to \cite{CJZ}, \cite{JL} and \cite{KLP}, where the $\ff^\tau$-semimartingale decomposition of an $\ff$-martingale is recovered from it's $\gg^\tau$-semimartingale decomposition, we retrieve the $\gg$-semimartingale decomposition of an $\ff$-martingale from its $\wh \ff$-semimartingale decomposition.
However, the present work is different to the case of initial
and progressive enlargement on two levels. First, on the level of
assumptions:  to show that the hypothesis $(H')$ holds between
$\ff$ and $\gg$, we  assume  that, for $i = 1,\dots,k$, the
hypothesis $(H')$ holds between $\ff$ and $\ff^i$ instead of
assuming directly that the hypothesis $(H')$ holds between $\ff$
and $\wh \ff$ (In Jeanblanc and Le Cam \cite{JL},
Callegaro et al. \cite{CJZ} and Kchia et al. \cite{KLP}, the largest filtration $\wh \ff$ is the initially enlarged
filtration $\gg^\tau$). Secondly, on the level of computation: in our setting,
we cannot exploit the same techniques as the ones used in the case
of initial and progressive enlargement with a single random time. To be more specific, suppose we know that the $\gg^\tau$-semimartingale decomposition of an 
$\ff$-martingale $M$ is given by $Y + A$, where $Y$ is (for ease of demonstration) a $\gg^\tau$-martingale and $A$ is a $\gg^\tau$-predictable process of finite variation. Then for all locally bounded $\ff^\tau$-predictable process $V$,
\bde
\expect{}{V\centerdot M} = \expect{}{V\I_{\llb 0,\tau \rrb} \centerdot M} + \expect{}{V\I_{\rrb \tau,\infty \llb} \centerdot M}.
\ede
The first term can be computed by the Jeulin-Yor formula, while for the second term, we have by assumption
\bde
\expect{}{V\I_{\rrb \tau,\infty \llb} \centerdot M} = \expect{}{V\I_{\rrb \tau,\infty \llb} \centerdot A}.
\ede
Formally, in order to make sure that the finite variation part is $\ff^\tau$-adapted, one should perform one more step in the calculation and take the $\ff^\tau$-dual projection of $A$. However, the computation stops here, since in the case of initial and progressive enlargement, one can exploit the fact that the process $\I_{\rrb \tau,\infty \llb} \centerdot A$ is $\ff^\tau$-predictable (see Lemma 2. in \cite{KLP}).

In our setting, suppose that the $\wh \ff$-predictable process of finite
variation in the $\wh\ff$-semimartingale decomposition of the $\ff$-martingale $M$ is known, which we again denoted by $A$. The differences here is that $\I_{D_i} \centerdot A$ is not necessarily $\gg$-predictable and theoretically, in order to find the $\gg$-semimartingale decomposition of $M$, one is forced to compute the $\gg$-dual predictable
projection of $A$. However in this paper, we work with the $\gg$-optional projection.
It is technically difficult and one soon realizes that the
computation is not trivial and a large part of this subsection is
devoted to overcome the technical difficulties associated with
localization and integrability conditions  that  one encounters
when computing the projections.

\subsubsection{Stopping times and Increasing processes: general results}
Before proceeding, we point out that for the purpose of this paper,
increasing processes are positive and finite valued unless
specified otherwise. In the following, we give two technical
results (\cref{t3} and \cref{vv}) for two given filtrations
$\kk$ and $\hh$ satisfying the usual conditions and a subset $D$ of $\Omega$ such that \be
\label{kh}\mathcal{P}(\kk)\cap D \subset \mathcal{P}(\hh)\cap
D\,.\ee 
We show how one can associate with every $\kk$-stopping
time (resp. $\kk$-increasing process) an $\hh$-stopping time (resp.
$\hh$-increasing process which can take value increasing) such
that they are equal on the set $D$.

In \cref{t3} and \cref{vv}, we set
\begin{align}
\wt N &= \,^{o,\hh}(\I_D) ,\quad \quad R_n  = \inf\seq{t\geq 0: \wt N_t \leq \frac{1}{n}}
\end{align}
for $n \in \nn^+$ and $R := \sup_n R_n = \inf\seq{s:\wt N_{s}\wt
N_{s-}= 0}$.

\bl \label[lem]{t3} 
Assume that $\hh$ and $\kk$ satisfy (\ref{kh}).
For any increasing sequence of $\kk$-stopping times
$(T_n)_{n\in \nn^+}$, there exists an increasing sequence of
$\hh$-stopping times $(S_n)_{n\in \nn^+}$ such that $T_n\I_D =
S_n\I_D$. In addition, if $\sup_nT_n=\infty$, then \hfill\break
(i) $S := \sup_n S_n$ is greater or equal to $R$.\hfill\break
 (ii) $\cup_n\seq{R_n = R} \subset \cup_n\seq{S_n \geq R}$.
\el
\begin{proof}
 Let $(T_n)_{n\in \nn^+}$ be an increasing sequence of $\kk$-stopping times. For every $n$, there exists an $\hh$-predictable process $H_n$ such that $\I_{(T_n,\infty)}\I_D=H_n\I_D$. Replacing $H_n$ by $H_n^+\wedge 1$, we can suppose that $0\leq H_n\leq 1$. Replacing $H_n$ by $\prod_{k=1}^nH_k$, we can suppose that $H_n\geq H_{n+1}$. Let $S_n=\inf\{t\geq 0: H_n(t)=1\}$.
 Then,  $(S_n)_{n\in \nn}$ is an increasing sequence of $\hh$-stopping times. We note that $S_n=T_n$ on the set $D$. \hfill\newline
(i) If $\sup_nT_n=\infty$,  by taking the $\hh$-predictable
projection, we see that
\bde
{^{p, \hh}}(\I_{\rrb T_n,\infty\llb
}\I_D)=\I_{\rrb S_n,\infty \llb}\wt {N}_-
\ede
then by using monotone convergence theorem and section theorem to the left-hand side, we obtain
\bde
 0=\I_{\cap_n
\rrb S_n,\infty \llb}\wt N_-.
\ede
For arbitrary $\epsilon >0$, we have from the above that $\wt N_{(S+\epsilon)-} = 0$, where $S = \sup_n S_n$. The process $\wt N$ is a positive $\hh$-supermartingale, this implies $S+\epsilon \geq R$ and therefore $S \geq R$, since $\epsilon$ is arbitrary.

(ii) Suppose there exists some $k$ such that $R_k = R$, then $\wt N_{R-} > 0$ and from $\I_{\cap_n \seq{S_n < R}}\wt N_{R-}=0$, one can there deduce that there exists $j$ such that $S_j \geq R$.
\end{proof}

The goal in the following is to study the measures associated with
finite variation process considered in different filtrations. As
in Jacod \cite{Ja}, we   use the concept of dominated measures, in
order to define a measure associated with a process which is the
difference of two increasing unbounded processes.

\bl\label[lem]{vv}
There exists a map $\wh V \longrightarrow \psi(\wh V)$ from the
space of $\kk$-predictable increasing processes to the space of
$\hh$-predictable increasing processes, valued in   $\rr_+ \cup
\{+\infty\}$, such that the following properties hold \hfill\break
(i) $\I_D \wh V=\I_D \psi(\wh V)$ and the support of the measure
$d\psi(\wh V)$ is contained in $\cup_{n}{\llb 0,
R_n\rrb}$\hfill\break
(ii) for any $\kk$-stopping time $\wh U$, there exists an
$\hh$-stopping time $U$ such that \bde \psi(\I_{\llb 0,\wh
U\rrb}\centerdot\wh V) = \I_{\llb 0, U\rrb}\centerdot\psi(\wh V)
\ede (iii) if the process $\wh V$ is bounded, then $\psi(\wh V)$
is  bounded by the same constant.
\el

\begin{proof}
Let $\wh V$ be an increasing $\kk$-predictable process. By
\cref{pp}, there exists an $\hh$-predictable process $V$ such that
$\wh V\I_D=V\I_D$.

(i) Given a $\kk$-localizing sequence $(T_n)_{n\in \nn^+}$ such
that $\wh V^{T_n}$ is bounded, by \cref{t3} there exists an
increasing sequence of $\hh$-stopping times $(S_n)_{n\in \nn^+}$
such that $\wh V^{T_n}\I_D = V^{S_n}\I_D$. By taking the
$\hh$-optional projection, we have \bde ^{o,\hh}(\wh V^{T_n}\I_D)
= \,^{o,\hh}(V^{S_n}\I_D)  = V^{S_n}\wt N \,.\ede From Theorem 47,
in Dellacherie and Meyer \cite{DM2},   the process $^{o,\hh}(\wh
V^{T_n}\I_D)$ is c\`adl\`ag, which implies that for all $n, k \in
\nn^+$, the process $V$ is c\`adl\`ag on $\llb 0, S_n\rrb\cap\llb
0, R_k\rrb$. Therefore $V$ is c\`adl\`ag on $\cup_k\llb 0,
R_k\rrb$, since, due to the fact that $\sup_n S_n \geq R$, for all
$k$, $\llb 0, R_k\rrb\subset \cup_n\llb 0, S_n\rrb$.

For every $n \in\nn^+$ and any pair of rationals $s \leq t \leq
R_n$, we have $\P$-a.s \bde V_t\I_{D} = \wh V_{t}\I_{D} \geq \wh
V_{s}\I_{D} = V_{s}\I_{D} \ede and by taking the conditional
expectation with respect to $\H_{t-}$, we have $\P$-a.s, the
inequality $V_{t}\wt N_{t-} \geq V_s\wt N_{t-}$. Using the fact
that the process $V$ is c\`adl\`ag on $\llb 0,R_n\rrb$, we can
first take right-limit in $t$ to show that for all $t \in
\rr_+\cap {\llb 0,R_n\rrb}$ and $s \in \Q_+\cap {\llb 0,R_n\rrb}$
the inequality $V_{t} \geq V_s$ holds. Then by taking limit in $s$
we extend this inequality to all $s\leq t \leq R_n$.

From the above, one deduces that the process $V$ is c\`adl\`ag and
increasing on $\cup_n\llb 0,R_n\rrb$. Then, one  defines the
increasing process $\psi(\wh V)$, which may take the value
infinity by setting \be \label{psiv} \psi(\wh V):=
{V}\I_{\cup_n\llb 0,R_n\rrb}+\lim_{s\uparrow R}{V}_s\I_{\{\wt
{N}_{R-}=0, 0<R<\infty\}}\I_{\llb R,\infty\llb}+{V}_R\I_{\rrb
R,\infty\llb}\I_{\{\wt {N}_{R-}>0, 0<R<\infty\}}, \ee which is
$\hh$-predictable since $\{\wt {N}_{R-}=0, 0<R<\infty\}\cap {\llb
R,\infty\llb}$ is the intersection of the set $(\cup_n \llb 0,R_n
\rrb)^c$ and the complement of ${\rrb R,\infty\llb}\cap{\{\wt
{N}_{R-}>0, 0<R<\infty\}}$, which are both $\hh$-predictable. From
\cref{psiv}, we see that the support of $d\psi(\wh V)$ is
contained in $\cup_n\llb 0,R_n\rrb$ and on the set $D$, we have
$\sup_n R_n= R = \infty$ which implies $\wh V\I_D = \psi(\wh
V)\I_D$.

(ii) For any $\kk$-stopping time $\wh U$ and any $\kk$-predictable
increasing process $\wh V$, we have from (i), the equality
$\psi(\I_{\llb 0, \wh U\rrb}\centerdot\wh V)\I_D= (\I_{\llb 0, \wh
U\rrb}\centerdot\wh V)\I_D$ and
\begin{align*}
(\I_{\llb 0, \wh U\rrb}\centerdot\wh V)\I_D= \I_{\llb 0,U\rrb}\centerdot(\wh V\I_D)&= \I_{\llb 0, U\rrb}\centerdot(\psi(\wh V)\I_D)=(\I_{\llb 0, U\rrb}\centerdot\psi(\wh V))\I_D
\end{align*}
where the existence of the $\hh$-stopping time $U$ in the first
equality follows from \cref{t3} and the second equality follows
from (i). By taking the $\hh$-predictable projection, we conclude
that for every $n \in \nn^+$, the processes $\I_{\llb
0,U\rrb}\centerdot\psi(\wh V)$ and $\psi(\I_{\llb 0,\wh
U\rrb}\centerdot \wh V)$ are equal on $\llb 0,R_n\rrb$ and
therefore on $\cup_n\llb 0,R_n\rrb$. This implies that the
processes $\I_{\llb 0,U\rrb}\centerdot\psi(\wh V)$ and
$\psi(\I_{\llb 0,\wh U\rrb}\centerdot \wh V)$ are equal
everywhere, since they do not increase on the complement of
$\cup_n\llb 0,R_n\rrb$.

(iii) If the process $\wh V$ is bounded by $C$, then for every $n
\in \nn^+$, we have $V \leq
\frac{\,^{p,\gg}(C\I_D)}{\,^{p,\gg}(\I_D)} = C$ on every interval
$\llb 0, R_n \rrb$ and therefore on $\cup_n\llb 0, R_n \rrb$. We
deduce from\cref{psiv} that on the set $\cup_n\llb 0, R_n \rrb$
the process $\psi(\wh V)$ is bounded by $C$, which implies that
$\psi(\wh V)$   is bounded by $C$, since the support of $d\psi(\wh
V)$ is contained in the set $\cup_n\llb 0, R_n \rrb$.
\end{proof}

\brem\label[rem]{sign} Our original goal in \cref{vv} was to
define, starting from a $\kk$-predictable process of finite
variation $V = V_+ - V_-$, an $\hh$-predictable process  of finite
variation   by setting $\psi(\wh V) := \psi(\wh V_+) - \psi(\wh
V_-)$. However, this is problematic since from\cref{psiv} we see
that in general the processes $\psi(\wh V_\pm)$ can take the value
infinity at the same time. Therefore one can not make use of the
usual definitions. \erem

Following the assumption and notation of \cref{vv}, one can
associate  with a given  $\kk$-predictable process $\wh V$ of
finite variation ($\wh V=\wh V_+ - \wh V_-$),  a pair of
$\hh$-predictable increasing processes $\psi(\wh V_+)$ and
$\psi(\wh V_-)$. In order to treat the problem mentioned in the
above remark, we define an auxiliary finite random measure $m$ on
$\Omega\times \mathcal{B}(\rr_+)$ by setting \be\label{commonm} 
dm
:=  (1 + \psi(\wh V_+) + \psi(\wh V_-))^{-2} d(\psi(\wh
V_+)+ \psi(\wh V_-)). \ee Since the processes $\psi(\wh V_+)$ and
$\psi(\wh V_-)$ can only take value infinity at the same time, we
deduce that $d\psi(\wh V_{\pm})$ is absolutely continuous with
respect to $dm$, which is absolutely continuous with respect to
$d(\psi(\wh V_+)+ \psi(\wh V_-))$. By (i) of \cref{vv},  this
implies that the support of $m$ is contained in $\cup_n \llb 0,R_n
\rrb$. Let us denote by $q^{\pm}$ the Radon-Nikod\'ym density of
$d\psi(\wh V_{\pm})$ with respect to $m$ and introduce the
following space of $\Omega\times \mathcal{B}(\rr_+)$-measurable
functions, 
\be\label{psivv1} \seq{f :\forall t > 0, \,\int_{(0,t]} |f_s||q^{+}_s -
q^{-}_s|\, m(ds) < \infty  }. \ee
 We define a linear operator
$\psi(\wh V)$ on the above set (\ref{psivv1}), which maps $f $ to
a process by setting for every $t \geq 0$, \be\label{psivv} f\ast
\psi(\wh V)_t: = \int_{(0,t]} f_s(q^{+}_s - q^{-}_s) \, m(ds). \ee
We shall call $\mathcal{L}^1(\psi(\wh V))$ the set (\ref{psivv1})
and an $\Omega\times \mathcal{B}(\rr)$-measurable function $f$ is
said to be $\psi(\wh V)$ integrable if $f\in
\mathcal{L}^1(\psi(\wh V))$. One should point out that the measure
$m$ is introduced to avoid the problem mentioned in \cref{sign}
and  that the set defined in (\ref{psivv1}) and the map defined in
(\ref{psivv}) are essentially independent of the choice of $m$.


\subsubsection{The $\gg$-Semimartingale Decomposition}
In this subsection, we place ourselves in the setting of \cref{sH12}, and we are in
position to derive the $\gg$-semimartingale decomposition of $\ff$-martingales. Let us first summarize the previous
results and introduce new notations.

We suppose that for each $i = 1,\dots, k$, the
$\ff^i$-semimartingale decomposition of an $\ff$-martingale $M$ is
given by $M = M^i + K^i$, where $M^i$ is an $\ff^i$-local
martingale and $K^i$ an $\ff^i$-predictable process of finite
variation. Then, from \cref{decompminmax}, for every $i = 1,\dots
,k$, the process
\begin{align}\label{hatm}
 \wh M^i     &:= M\I_{D_i} - \I_{D_i}\left(K^i\I_{D_i} +
\frac{\I_{D_i}}{N^i_-}\centerdot\pbracket{N^i}{M}^i\right)
\end{align}
is an $\wh \ff$-local martingale. For simplicity, for every $i= 1,\dots, k$, we adopt the following notation:
\begin{align}
\wt N^i &:= \,^{o,\gg}(\I_{D_i}), \quad N^i := \,^{o,\ff^i}(\I_{D_i}) \label{tilden}\\
\wh V^i &:= K^i+
\frac{1}{N^i_-}\centerdot\pbracket{N^i}{M}^i.\label{hatv}
 \end{align}
From \cref{pp}, for every $i =1,\dots,k$, we have $\mathcal{P}(\wh
\ff)\cap D_i \subset \mathcal{P}(\gg)\cap D_i$ and one can   apply
\cref{vv} with the set $D_i$, the filtrations $\wh \ff$ and $\gg$
and define the linear operator $\psi(\wh V^i)$ as in (\ref{psivv})
on the space of  $\psi(\wh V^i)$ integrable functions given by
\bde \mathcal{L}^1(\psi(\wh V^i)) := \seq{f:
\int_{\rr_+}|f_s||q^{i,+}_s - q^{i,-}_s|dm^i_s <\infty }, \ede
where the measure $m^i$ is constructed from the $\wh
\ff$-predictable process $\wh V^i$ as shown in (\ref{commonm}) and
$q^{i,\pm}$ is the density of $\psi(\wh V^i_{\pm})$ with respect
to $m^i$.

For an arbitrary filtration $\kk$, we write $X\stackrel{\kk\rm
-mart} = Y$, if $X-Y$ is a $\kk$-local martingale.

\bt\label{gde} Let  $M$ be  $\ff$-martingale, such that for every
$i= 1,\dots k$, $M= M^i + K^i$, where $M^i$ is an $\ff^i$-local
martingale and $K^i$ an $\ff^i$-predictable process of finite
variation. Then under \cref{g}, for every $i = 1,\dots, k$, the
process $\wt N^i$ defined in (\ref{tilden}) belongs to
$\mathcal{L}^1(\psi(\wh V^i))$, where $\wh V^i$ is defined in
(\ref{hatv}) and
\begin{align}\label{dc1}
M -  \sum^k_{i=1} \wt N^i_-\ast \psi(\wh V^i)
\end{align}
is a $\gg$-local martingale, where $\ast$ is defined in
(\ref{psivv}). \et
\begin{proof}
It is sufficient to show that for any fixed $i = 1,\dots ,k$, the $\gg$-martingale $\wt N^i$ is in $\mathcal{L}^1(\psi(\wh V^i))$ and $M\wt N^i - \wt N^i_-\ast \psi(\wh V^i)$ is an $\gg$-local martingale.

Before proceeding, we point out that the process $M$ is by \cref{fhatdecomp} an $\wh \ff$-semimartingale 
and therefore an $\gg$-semimartingale by the result of Stricker \cite{St}. To see that $M$ is an $\gg$-special semimartingale, 
we note that there exists a sequence of $\ff$-stopping times (therefore $\gg$-stopping times) 
which reduces $M$ to $\mathcal{H}^1$, and one can check directly, that $M$ is an $\gg$-special 
semimartingale by using special semimartingale criteria such as Theorem 8.6 in \cite{HWY}. 
The aim in the rest of the proof is to find explicitly the $\gg$-semimartingale decomposition of $M$.

The process $M\wt N^i$ is an $\gg$-special semimartingale
since $M$ is an $\gg$-special semimartingale and $\wt N^i$ is a
bounded $\gg$-martingale. Let us denote $B^i$ the unique
$\gg$-predictable process of finite variation in the
$\gg$-semimartingale decomposition of $M\wt N^i$.

For every $n \in \nn^+$, we define $R_n = \inf\seq{t\geq 0, \wt
N^i \leq \frac{1}{n}}$ and $R =\sup_n R_n = \inf\seq{t\geq 0, \wt
N^i_-\wt N^i =0}$. The method of the proof is to identify the
process $B^i$ with $-\wt N^i_-\ast\psi(\wh V^i)$. To do that, it is
sufficient to show that the two processes coincide on the sets
$\cup_n\llb 0, R_n\rrb$ and $\llb 0, R\rrb\setminus \cup_n\llb 0,
R_n\rrb$: indeed  $B^i$ and $\wt N^i_-\ast \psi(\wh V^i)$ are
constant on $\rrb R, \infty \llb$, as $M_t\wt N^i_t = 0$ for $t>
R$ and the support of $dm^i$ is in the complement of $\rrb R,
\infty \llb$ by (i) of \cref{vv}.

On the set $\cup_n\llb 0, R_n\rrb$, let $(T_n)_{n \in \nn^+}$ be a
localizing sequence of $\wh \ff$-stopping times such that the
process $\wh M^i$ defined in (\ref{hatm}) stopped at $T_n$, i.e.,
$(\wh M^i)^{T_n}$ is an $\wh \ff$-martingale and $(\wh V^i)^{T_n}$
is of bounded total variation, then
\begin{align}\label{gg}
(\wh M^{i})^{T_n} & = M^{T_n}\I_{D_i} + \I_{\llb 0, T_n\rrb}\centerdot\wh V^{i}_+\I_{D_i} - \I_{\llb 0, T_n\rrb}\centerdot\wh V^{i}_-\I_{D_i}\\
            & = M^{S_n}\I_{D_i} + \I_{\llb 0, S_n\rrb}\centerdot\psi(\wh V^{i}_+)\I_{D_i} - \I_{\llb 0, S_n\rrb}\centerdot\psi(\wh V^{i}_-)\I_{D_i}\nonumber
\end{align}
where the second equality holds from (i) of \cref{vv} with the
existence of the sequence of $\gg$-stopping times $(S_n)_{n\in
\nn^+}$ given by \cref{t3}. Since for every $n \in \nn^+$, the
$\wh \ff$-adapted processes $\I_{\llb 0, T_n\rrb}\centerdot\wh
V^{i}_\pm$ are bounded, by (ii) and (iii) of \cref{vv}, we can
conclude that the processes $\I_{\llb 0, S_{n}
\rrb}\centerdot\psi(\wh V^i_\pm)$ are bounded $\gg$-predictable
processes. Together with the property that the Lebesgue integral
and the stochastic integral coincide when all quantities are
finite, we obtain that 
\bde \I_{\llb 0, S_n\rrb}\centerdot\psi(\wh
V^{i}_+) - \I_{\llb 0, S_n\rrb}\centerdot\psi(\wh V^{i}_-) =
\I_{\llb 0, S_n\rrb}\ast \psi(\wh V^{i}), \ede
 that is the
left-hand side coincides with the application of $\psi(\wh V^{i})$
on $\I_{\llb 0, S_n\rrb}$ defined in\cref{psivv}. By taking the
$\gg$-optional projection, we obtain
\begin{align*}
M^{S_n}\wt N^i    \stackrel{\gg\rm -mart} = - \wt N^i\left(
\I_{\llb 0, S_{n} \rrb}\ast\psi(\wh V^i)\right) & \stackrel{\gg\rm -
mart}=  - \wt N^i_-\centerdot\left(\I_{\llb 0, S_{n}
\rrb}\ast\psi(\wh V^i)\right)
\end{align*}
where the second equality follows from integration by parts and
Y\oe urp's lemma. For each $n \in \nn^+$, from uniqueness of the
$\gg$-semimartingale decomposition, we deduce that  
\bde \I_{\llb
0, S_n \rrb} \centerdot B^i = - \wt
N^i_-\centerdot\left(\I_{\llb 0, S_{n} \rrb}\ast\psi(\wh
V^i)\right)= - (\wt N^i_- \I_{\llb 0, S_{n} \rrb} )\ast \psi(\wh
V^i), \ede where the second equality follows from the property
that the Lebesgue integral and the stochastic integral coincide
when all quantities are finite. Then for any fixed $t\geq 0$ and
$H_s := \I_\seq{q^{i,+}_s - q^{i,-}_s < 0}- \I_\seq{q^{i,+}_s -
q^{i,-}_s \geq 0}$, \bde \int_{(0,t]} H_s\I_\seq{s\leq S_n}
dB^i_s = \int_{(0,t]}  \I_\seq{s \leq S_{n}} \wt N^i_-|q^{i,+}_s -
q^{i,-}_s|\,dm^i_s. \ede To take limit in $n$, one applies the
dominated convergence theorem to the left-hand side and the
Beppo-Levi (monotone convergence) theorem to the right. One then
concludes that $\wt N^i\I_{\cup_n \llb 0,S_n \rrb}$ is $\psi(\wh
V^i)$ integrable as the limit on the left-hand side is finite.
This implies that $\wt N^i$ is $\psi(\wh V^i)$ integrable as the
support of $dm^i$ is contained in ${\cup_n \llb 0,R_n \rrb}$ which
is contained in ${\cup_n \llb 0,S_n \rrb}$ by (i) of \cref{t3}.

We have shown that $\wt N^i\in \mathcal{L}^1(\psi(\wh V^i))$ and therefore $\wt N^i_-\ast \psi(\wh V^i)$ is a process such that for all $n\in \nn^+$, we have $(\wt N^i_-\ast \psi(\wh V^i))^{S_n} = (B^i)^{S_n}$. This implies $\wt N^i_-\ast \psi(\wh V^i) = B^i$ on $\cup_n \llb 0,S_n\rrb$ and therefore on ${\cup_n \llb 0,R_n \rrb}$.

On the set $\llb 0, R\rrb\setminus \cup_n\llb 0, R_n\rrb$, one
only needs to pay attention to   the set $F  = \seq{\forall n, R_n
< R}= \seq{\wt N^i_{R-} = 0}$, since on the complement $F^c$, the
set $\llb 0, R\rrb\setminus \cup_n\llb 0, R_n\rrb$ is empty. From
Lemma 3.29 in He et al. \cite{HWY}, $R_{F} = R\I_F +
\infty\I_{F^c}$ is a $\gg$-predictable stopping time and on $F$,
we have $\llb 0, R\rrb\setminus \cup_n\llb 0, R_n\rrb = \llb
R_F\rrb$. From the fact that $\Delta(M\wt N^i)_{R_{F}} = 0$, Lemma
8.8 of \cite{HWY} shows that $|\Delta B^i_{R_{F}}|= 0$. On the
other hand, from \cref{psiv}, we deduce that the measure $dm^i$,
which is absolutely continuous with respect to $\psi(\wh V^i_+) +
\psi(\wh V^i_-)$ has no mass at $\llb R_F\rrb$ and therefore $\wt
N^i_-\ast\psi(\wh  V^i)$ does not jump at $\llb R_F \rrb$. This
implies that the jumps of the processes $B^i$ and $\wt
N^i_-\ast\psi(\wh  V^i)$ coincide on $\llb 0, R\rrb\setminus
\cup_n\llb 0, R_n\rrb$.
\end{proof}

We conclude this section with the following lemma, which appears to be very useful when one wants to compute the $\gg$-drift in practice. 

\bl\label[lem]{gdel}
For $1 \leq i \leq k$, the process $\wt N^i_-\ast \psi(\wh V^i)$ is the dual $\gg$-predictable projection of $\I_{D_i}\wh V^i$ and $N^i_-\centerdot \wh V^i$ is the $\ff^i$ drift of the $\ff^i$-special semimartingale $N^iM$. 
\el
\begin{proof}
Similar to the proof of \cref{gde}, for a fixed $i =1,\dots ,k$ and every $n \in \nn^+$, we define $R_n = \inf\seq{t\geq 0, \wt
N^i \leq \frac{1}{n}}$, $R =\sup_n R_n = \inf\seq{t\geq 0, \wt
N^i_-\wt N^i =0}$ and let $(T_n)_{n \in \nn^+}$ be a
localizing sequence of $\wh \ff$-stopping times (and $(S_n)_{n\in \nn^+}$ the corresponding $\gg$-stopping times from \cref{t3}) such that the
process $\wh M^i$ defined in (\ref{hatm}) stopped at $T_n$, i.e.,
$(\wh M^i)^{T_n}$ is an $\wh \ff$-martingale, and $(\wh V^i)^{T_n}$
is of bounded total variation. Again, the method of the proof is to identify the process $(\wh V^{i}\I_{D_i})\,^{p,\gg}$ with $-\wt N^i_-\ast\psi(\wh V^i)$.

By taking the $\gg$-optional projection in \eqref{gg}, we see that for a fixed $i = 1,\dots, k$, 
\begin{align*}
M^{S_n}\wt N^i  \stackrel{\gg\rm -mart} = - \,^{o,\gg}(\I_{\llb 0, T_n\rrb}\centerdot \wh V^{i}\I_{D_i}) & \stackrel{\gg\rm -mart} = - (\I_{\llb 0, T_n\rrb}\centerdot V^{i}\I_{D_i})\,^{p,\gg}\\
&\quad = - (\I_{\llb 0, S_n\rrb}\centerdot \wh V^{i}\I_{D_i})\,^{p,\gg}\\
													 &\quad = - \I_{\llb 0, S_n\rrb}\centerdot (\wh V^{i}\I_{D_i})\,^{p,\gg},
\end{align*}
where the second and last equality follows from Corollary 5.31 and Corollary 5.24 in \cite{HWY} respectively. From \cref{gde} and the uniqueness of $\gg$-special semimartingale decomposition, the process $\wt N^i_-\ast \psi(V)$ is equal to $(\wh V^{i}\I_{D_i})\,^{p,\gg}$ on $\llb 0,S_n\rrb$ for all $n$ and therefore on $\cup_{n}\llb 0, R_n \rrb$. On the complement $(\cup_{n}\llb 0, R_n \rrb)^c$, for every $C \in \F$ and every bounded $\gg$-predictable process $\xi$, by duality
\bde
\expect{\P}{\I_C(\xi\I_{(\cup_{n}\llb 0, R_n \rrb)^c}\centerdot (\wh V^{i}\I_{D_i})\,^{p,\gg})}= \expect{\P}{(\,^{p,\gg}(\I_C)\xi\I_{\cap_{n}\rrb R_n, \infty \llb}\centerdot \wh V^{i}\I_{D_i})} = 0 
\ede
where the last equality holds, since on the set $D_i$, $\sup_{n}R_n = \infty$. This implies $(\wh V^{i}\I_{D_i})\,^{p,\gg}$ also coincide with $\wt N^i_-\ast \psi(V)$ on the $(\cup_{n}\llb 0, R_n \rrb)^c$ and therefore everywhere.\hfill\break\hfill\break
To show for every $i=1,\dots, k$, the process $N^i_-\centerdot \wh V^i$ is the $\ff^i$-drift of the $\ff^i$-special semimartingale $N^iM$, it is sufficient to apply the $\ff^i$-integration by parts formula to $N^iM = N^i(M-K^i) + N^iK^i$. We see that
\begin{align*}
N^iM &\stackrel{\ff^i\rm -mart}= \pbracket{N^i}{M-K^i}^i + N^i_-\centerdot K^i	       \stackrel{\ff^i \rm-mart}= \pbracket{N^i}{M}^i + N^i_-\centerdot K^i,
\end{align*}
where the second equality follows from Y\oe urp's lemma. One can now conclude that the $\ff^i$-drift of $N^iM$ is given by $N^i_-\centerdot \wh V^i$ from the uniqueness of the $\ff^i$-special semimartingale decomposition.
\end{proof}

\vfill\break







\

\section{Application to Multiple Random Times}\label{applied}
Our goal here is to apply the previous methodology to enlargement
of filtrations with random times.

\subsection{Progressive enlargement with random times and their re-ordering}\label{multi}

We introduce the following notations. For two elements $a,b$ in $[0,\infty]$ we denote $$
a\nmid b=
\left\{
\dcb
a&&\mbox{ if $a\leq b$,}\\
\\
\infty&&\mbox{ if $a> b$.}\\
\dce
\right.
$$
Let $\boldsymbol\xi=(\xi_1,\ldots,\xi_k)$ be an $k\in\mathbb{N}^*$ dimensional vector of random times. For $b\in[0,\infty]$, we write $\boldsymbol \xi\nmid b=(\xi_1\nmid b,\ldots,\xi_k\nmid b)$ and $
\sigma(\boldsymbol\xi \nmid b) =
\sigma(\xi_1\nmid b)\vee\dots \vee \sigma(\xi_k\nmid b).
$

\bd The progressive enlargement of $\ff$ with the family of random
times $\boldsymbol\xi$ is denoted
$\ff^{\boldsymbol\xi} = (\F^{\boldsymbol\xi}_t)_{t\geq
0}$; this is the smallest filtration satisfying the usual conditions containing $\ff$ making $\xi_1,\dots,\xi_k$ stopping times. In other terms, \bde 
\F^{\boldsymbol\xi}_t :=
\F^{\xi_1,\dots,\xi_k}_t :=
\bigcap_{s>t}\left(\F_s\vee
\sigma(\boldsymbol\xi \nmid s)\right),\ t\in\mathbb{R}_+. 
\ede 
\ed

\bd The initial enlargement of $\ff$ with the family of random times
$\boldsymbol\xi$ is denoted $\ff^{\sigma(\boldsymbol\xi)}:=\gg^{\xi_1,\dots,\xi_k} =
(\F^{\sigma(\boldsymbol\xi)}_t)_{t\geq 0}$; this  is the smallest
filtration containing $\ff$, satisfying the usual conditions,
such that the random times $\boldsymbol\xi=(\xi_1,\dots,\xi_k)$ are
$\F^{\sigma(\boldsymbol\xi)}_0$ measurable. One has \bde
\F^{\sigma(\boldsymbol\xi)}_t=
\bigcap_{s>t}\left(\F_s\vee\sigma(\boldsymbol\xi)\right),\ t\in\mathbb{R}_+. 
\ede 
\ed

We need a sorting rule (cf. \cite{S}). For any function $\mathfrak{a}$ defined on $\{1,\ldots,k\}$ taking values $\{a_1,\ldots,a_k\}$ in $[0,\infty]$, consider the points $(a_1,1),\ldots,(a_k,k)$ in the space $[0,\infty]\times\{1,\ldots,k\}$. These points are two-by-two distinct. We order these points according to the alphabetic order in the space $[0,\infty]\times\{1,\ldots,k\}$. Then, for $1\leq i\leq k,$ the rang of $(a_i,i)$ in this ordering is given by $$
R^\mathfrak{a}(i)=R^{\{a_1,\ldots,a_k\}}(i)=\sum_{j=1}^k \ind_{\{a_j<a_i\}}+\sum_{j=1}^k \ind_{\{j<i,a_j=a_i\}}+1.
$$
The map $i\in\{1,\ldots,k\}\rightarrow R^\mathfrak{a}(i)\in\{1,\ldots,k\}$ is a bijection. Let $\rho^\mathfrak{a}$ be its inverse. Define $\uparrow\!\!\!\mathfrak{a}=\mathfrak{a}(\rho^\mathfrak{a})$, where $\uparrow\!\!\!\mathfrak{a}(j)$ can be roughly qualified as the $j$th value in the increasing order of the values $\{a_1,\ldots,a_k\}$. $\uparrow\!\!\!\mathfrak{a}$ is an non decreasing function on $\{1,\ldots,k\}$ taking the same values as $\mathfrak{a}$.

We consider $\omega$ by $\omega$ the random function $\mathfrak{a}^{\boldsymbol\xi}$ on $\{1,\ldots,k\}$ taking the values $\{\xi_1,\ldots,\xi_k\}$ and define the non decreasing function $\uparrow\!\!\!\mathfrak{a}^{\boldsymbol\xi}$ as above.

\bl\label{re-ordering}
For any $1\leq j\leq k$, there exists a Borel function $\mathfrak{s}_{j}$ on $[0,\infty]^k$ such that $$
\uparrow\!\!\!\mathfrak{a}^{\boldsymbol\xi}=\mathfrak{s}_{j}(\xi_1,\ldots,\xi_k).
$$
If the $\xi_1,\ldots,\xi_k$ are stopping times with respect to some filtration, the random times $\uparrow\!\!\!\mathfrak{a}^{\boldsymbol\xi}(1),\ldots,\uparrow\!\!\!\mathfrak{a}^{\boldsymbol\xi}(k)$ also are stopping times with respect to the same filtration.
\el

\begin{proof} This is a consequence of the following identity\!\!: for any $t\geq 0$,$$
\{\uparrow\!\!\!\mathfrak{a}^{\boldsymbol\xi}(j)\leq t\}=\cup_{I\subset\{1,\ldots,k\}, \sharp I=j}\{\xi_h\leq t, \forall h\in I\}. 
$$
\end{proof}

We will call the random times $\uparrow\!\!\!\mathfrak{a}^{\boldsymbol\xi}(1),\ldots,\uparrow\!\!\!\mathfrak{a}^{\boldsymbol\xi}(k)$ the increasing re-ordering of $\{\xi_1,\ldots,\xi_k\}$, and we denote $\xi_{(j)}=\uparrow\!\!\!\mathfrak{a}^{\boldsymbol\xi}(j), 1\leq j\leq k$, and $\overline{\boldsymbol\xi}=(\xi_{(1)},\ldots,\xi_{(k)})$. We define then $\ff^{\overline{\boldsymbol\xi}}$ and $\ff^{\sigma(\overline{\boldsymbol\xi})}$. Also we define $\xi_{(0)}=0,\xi_{(k+1)}=\infty$ to complete the increasing re-ordering (only when the length of the initial random vector $\boldsymbol \xi$ is $k$).

The following result will be useful.

\bl\label{rtql}
Suppose $\btau$ and $\boldsymbol\xi$ are two $k$-dimensional
vectors of random times. If a measurable set $D$ is a subset of $\seq{\btau=\boldsymbol\xi}$,
then   $\F^\btau_t\cap D=\F^{\boldsymbol\xi}_t\cap D$ for every
$t\geq 0$. \el

\begin{proof}

By symmetry, it is enough to show that for all $t\geq 0$, we have
$\F^\btau_t\cap D \subset \F^{\boldsymbol\xi}_t\cap D$. Fix $t\geq
0$, then for any $A \in \F^\btau_t$ and all $n\geq 1$, there exist
an $\F_{t+1/n}\otimes \mathcal{B}(\rr^d)$-measurable function
$g_n$ such that $\I_A = g_n(\tau \wedge (t+\frac{1}{n}))$. We have \bde \I_A \I_D = g_n(\btau\wedge (t+1/n))\I_D =
g_n(\boldsymbol \xi\wedge (t+1/n))\I_D 
\ede 
and then $$
A\cap
D = \seq{\lim _{n\rightarrow \infty}g_n(\boldsymbol\xi\wedge
(t+1/n)) = 1}\cap D
$$
whilst $\seq{\lim _{n\rightarrow
\infty}g_n(\xi\wedge (t+1/n)) = 1}
$ 
is
$\F^{\boldsymbol\xi}_t$-measurable by right continuity of
$\ff^{\boldsymbol\xi}$.
\end{proof}

The following notations will be used. Let $n$ be an integer $n\geq k$. Let $I_k=\{1,\ldots,k\}$ and $I_n=\{1,\ldots,n\}$. For a subset $J\subset I_n$ denote $\boldsymbol\xi_J=(\xi_j:j\in J)$ in their natural order. In particular, if $J=I_k$ we denote $\boldsymbol\xi_{I_k}$ be $\boldsymbol\xi_k$. For any injective map $\varrho$ from $I_k$ into $I_n$ denote $\boldsymbol\xi_\varrho=(\xi_{\varrho(j)}:j\in I_k)$ and $\boldsymbol\xi_{\setminus\varrho}=(\xi_{j}:j\in I_n\setminus \varrho(I_k))$ in their natural order.

\

\subsection{Hypothesis $(H')$ for ordered random times}\label{orderedtimes}

From now on we are given an $n\in\mathbb{N}^*$ dimensional vector of random times $\boldsymbol\tau=(\tau_1,\ldots,\tau_n)$. Let $k$ be an integer $1\leq k\leq n$. The aim of this section is to apply the results from \cref{sH1}
to investigate the hypothesis
$(H')$ between $\ff$ and $\ff^{\overline{\boldsymbol\tau}_k}=\ff^{\tau_{(1)},\ldots,\tau_{(k)}}$.

Let $\mathcal{I}(k,n)$ denote the set of all injective functions $\varrho$ from $I_k$ into $I_n$. For $\varrho \in \mathcal{I}(k,n)$, we
introduce the set $d_{\varrho} :=
\seq{(\tau_{(1)},\dots,\tau_{(k)})= (\tau_{\varrho(1)}, \dots,
\tau_{\varrho(k)})}$. Let $\delta$ be a strictly increasing map from $\mathcal{I}(k,n)$ into $\mathbb{N}$. We define a  partition of
the space $\Omega$  by setting $$D_\varrho := d_\varrho\cap
(\cup_{\delta(\varrho')<\delta(\varrho)} d_{\varrho'})^c\,.
$$
Following Lemma \ref{rtql} we have

\bcor\label[cor]{cmult} The triplet $(\ff,
(\ff^{\btau_\varrho})_{\varrho \in \mathcal{I}(k,n)},
\ff^{\overline{\boldsymbol\tau}_k})$ satisfies \cref{g} with
respect to the partition $(D_\varrho)_{\varrho \in
\mathcal{I}(k,n)}$. 
\ecor

With the family $(\ff^{\btau_\varrho})_{\varrho \in \mathcal{I}(k,n)}$ and the partition $(D_\varrho)_{\varrho \in
\mathcal{I}(k,n)}$, we can construct the filtration $\wh
\ff$ as in \cref{cmult}, which will be denoted by $\wh
\ff^{\overline{\boldsymbol\tau}_k}=(\widehat\F^{\overline{\boldsymbol\tau}_k}_t)_{t\in\mathbb{R}_+}$, by setting for every
$t \geq 0$, 
\be 
\widehat\F^{\overline{\boldsymbol\tau}_k}_t :=
\seq{A\in \F\, |\,\forall \varrho \in \mathcal{I}(k,n), \,\exists
A^{\varrho}_t\in \F^{\boldsymbol\tau_\varrho}_t\,\,\mathrm{such \,\,that}\, A\cap D_{\varrho}
=  A^{{\varrho}}_t \cap D_{\varrho}}. 
\ee 
The $\wh \ff^{\overline{\boldsymbol\tau}_k}$
and $\ff^{\overline{\boldsymbol\tau}_k}$-semimartingale decompositions of
$\ff$-martingales are now readily available.

\bl\label[lem]{multiprop} Assume that for every $\varrho \in
\mathcal{I}(k,n)$, the $\ff$-martingale $M$ is an
$\ff^{\boldsymbol\tau_{\varrho}}$-semimartingale
with $\ff^{\boldsymbol\tau_{\varrho}}$-semimartingale decomposition given by $M =
M^\varrho+ K^\varrho$, where $M^\varrho$ is an
$\ff^{\boldsymbol\tau_{\varrho}}$-local
martingale and $K^\varrho$ is an $\ff^{\boldsymbol\tau_{\varrho}}$-predictable process of finite variation. We
denote
\begin{align*}
\wt N^\varrho &:= \,^{o,\ff^{\overline{\boldsymbol\tau}_k}}(\I_{D_\varrho}), \quad
N^\varrho := \,^{o,\ff^{\boldsymbol\tau_{\varrho}}}(\I_{D_\varrho})\\  
\wh V^\varrho & :=
K^\varrho +
\frac{1}{N^\varrho_-}\centerdot\pbracket{N^\varrho}{M}^\varrho,
\end{align*}
where the predictable bracket $\pbracket{.}{.}^\varrho$ is
computed with respect to the filtration
$\ff^{\boldsymbol\tau_{\varrho}}$. Then
\hfill\break (i) the $\wh \ff^{\overline{\boldsymbol\tau}_k}$-semimartingale
decomposition of $M$ is given by \bde M = \wh M + \sum_{\varrho
\in \mathcal{I}(k,n)}\I_{D_\varrho}\left(K^\varrho\I_{D_\varrho} +
\frac{\I_{D_\varrho}}{N^\varrho_-}\centerdot\pbracket{N^\varrho}{M}^\varrho\right)
\ede where  $\wh M$ is an $\wh \ff^{\overline{\boldsymbol\tau}_k}$-local martingale,
\hfill\break (ii) the $\ff^{\overline{\boldsymbol\tau}_k}$-semimartingale decomposition
of $M$ is given by
\begin{align*}
M = \wt M + \sum_{\varrho \in \mathcal{I}(k,n)} \wt N^\varrho_- \ast \psi(\wh V^{\varrho})
\end{align*}
where $\wt M$ is an $\ff^{\overline{\boldsymbol\tau}_k}$-local martingale and the linear operator $\psi(\wh V^{\varrho})$ is described in\cref{psivv}.
\el

\

\subsection{Computation under density hypothesis}\label{computation}
The formula in Lemma \ref{multiprop} allows one to compute the $\ff^{\overline{\boldsymbol\tau}_k}$-semimartingale decomposition of a $\ff$-martingale $M$. Compared with classical result on this subject, the formula in Lemma \ref{multiprop} has the particularity to express the drift of $M$ in $\ff^{\overline{\boldsymbol\tau}_k}$ as a {\it `weighted average'} of its $\ff^{\btau_{\varrho}}$-drifts when $\varrho$ runs over $\mathcal{I}(k,n)$. Such a weighted average formula may be useful for model analysis and risk management. The rest of this paper is devoted to illustrate this weighted average by an explicit computation in the case of density hypothesis. To this end we develop techniques which can have their utility elsewhere.

\bhyp\label[hyp]{density} The conditional laws of the vector $\btau = (\tau_1,\dots,
\tau_n)$ satisfies the density hypothesis. In other terms, for any non negative Borel function $h$ on $\rr^n_+$, for $t\in\mathbb{R}_+$, we have
$\P$-almost surely, 
\bde \cE{\P}{h(\btau)}{\F_t} =
\int_{\rr^n_+} h(\mathbf{u})a_t(\mathbf{u}) \mu^{\otimes
n}(\mathbf{du}), 
\ede 
where $\mu$ is a
non-atomic $\sigma$-finite measure on $\rr_+$, and $a_t(\omega,\mathbf{u})$ is a non negative $\mathcal{F}_t\otimes\mathcal{B}(\rr_+^n)$ measurable function, called the conditional density function at $t$. 
\ehyp

Remark that, according to \cite[Lemme(1.8)]{JJ}, the density function $a_t(\omega,\mathbf{x})$ can be chosen everywhere c\`adl\`ag in $t\in\mathbb{R}_+$, with however a rougher measurability : $a_t$ being $\cap_{s>t}(\mathcal{F}_s\otimes\mathcal{B}(\rr_+^n))$ measurable. Moreover, for fixed $\mathbf{x}$, $a_t(\mathbf{x}), t\in\rr_+,$ is a $(\mathbb{P},\mathbb{F})$ martingale. We assume this version of the density function in this section. Without loss of the generality, we assume $$
\int_{\rr^n_+} a_t(\mathbf{x}) \mu^{\otimes
n}(\mathbf{dx})=1,\ t\in\mathbb{R}_+,
$$
everywhere. Thus the regular conditional laws, denotes by $\nu_t$, of the vector $\btau$ with respect to the $\sigma$-algebra $\mathcal{F}_t$ has $a_t$ as the density function with respect to $\mu^{\otimes
n}$.

\bl\label{ahmu}
Let $t\in\rr_+$. For any bounded function $h$ on $\Omega\times\rr_+^n$, $\cap_{s>t}(\mathcal{F}_s\otimes\sigma(\rr_+^n))$ measurable, we have$$
\mathbb{E}[h(\btau)|\mathcal{F}_t]=\int_{\rr_+^n} a_t(\mathbf{x})h(\mathbf{x}) \mu^{\otimes n}(d\mathbf{x}).
$$
\el

\begin{proof}
For any $s>t$, $h$ is $\mathcal{F}_s\otimes\sigma(\rr_+^n)$ measurable. Let $B\in\mathcal{F}_t$.
$$
\begin{array}{lll}
\mathbb{E}[\ind_Bh(\btau)]
&=&\mathbb{E}[\int_{\rr_+^n}\ind_Bh(\mathbf{x})a_s(\mathbf{x})\mu^{\otimes n}(d\mathbf{x})]\\

&=&\mathbb{E}[\ind_B\int_{\rr_+^n}h(\mathbf{x})a_t(\mathbf{x})\mu^{\otimes n}(d\mathbf{x})],
\end{array}
$$
because, for every $\mathbf{x}$, $a_t(\mathbf{x}), t\in\rr_+,$ is a $(\mathbb{P},\mathbb{F})$ martingale, and $h(\mathbf{x})$ is $\mathcal{F}_t$ measurable, thanks to the right continuity of $\mathbb{F}$. The lemma is proved because $\int_{\rr_+^n}h(\mathbf{x})a_t(\mathbf{x})\mu^{\otimes n}(d\mathbf{x})$ also is $\mathcal{F}_t$ measurable.
\end{proof}

To make computations under \cref{density}, we introduce the following system of notations. For a subset $J\subset \{1,\ldots,n\}$ of cardinal $\#J=j\in\mathbb{N}$, for vectors
$\mathbf{z}\in\mathbb{R}^{j}$ and
$\mathbf{y}\in\mathbb{R}^{n-j}$, let $\mathfrak{c}_J(\mathbf{z},\mathbf{y})$ be the vector $\mathbf{x}\in\mathbb{R}^{n}$ for which $\mathbf{x}_J$ is given by $\mathbf{z}$ in their natural
order, and $\mathbf{x}_{J^c}$ is given by $\mathbf{y}$ in the natural order. For non negative Borel function $g$ we denote$$
g_J(\mathbf{z})=\int_{\rr^{n-j}_+} g(\mathfrak{c}_J(\mathbf{z},\mathbf{y}))\mu^{\otimes
(n-j)}(\mathbf{dy}).
$$
We check directly that $\btau_J$ satisfies the density hypothesis with density functions $(a_t)_J, t\in\rr_+$. We denote by $a_{J,t}, t\in\rr_+$, the c\`adl\`ag version of $(a_t)_J$ defined in \cite[Lemme(1.8)]{JJ}. Notice that, for all $t\in\rr_+$, for almost all $\omega$, $(a_t)_J(\omega,\mathbf{x})=a_{J,t}(\omega,\mathbf{x})$ $\mu^{\otimes n}$-almost everywhere. Hence, if $\mathfrak{p}_J$ denotes the projection map $\mathbf{x}\rightarrow\mathbf{x}_J$ on $\rr_+^n$, the conditional expectation of $g$ under $\nu_t(\omega)$ given $\mathfrak{p}_J=\mathbf{z}$ is the function $$
\mathbb{E}^{\nu_t}[g|\mathfrak{p}_J=\mathbf{z}]=\frac{(ga_t)_J(\mathbf{z})}{a_{J,t}(\mathbf{z})}\ind_{\{a_{J,t}(\mathbf{z})>0\}}.	
$$ 
For $\mathbf{x}\in\mathbb{R}^{n}$, let $\overline{\mathbf{x}}=(x_{(1)},\ldots,x_{(n)})$ be the increasing re-ordering of $\mathbf{x}$ and let $\mathfrak{r}$ be the map $\mathfrak{r}(\mathbf{x})=\overline{\mathbf{x}}$. Let $\mathfrak{S}$ be the symmetric group on $I_n$. For $\pi\in\mathfrak{S}$, we define the map $\pi(\mathbf{x})=(x_{\pi(1)},\ldots,x_{\pi(n)})$ and$$
\overline{\overline{g}}(\mathbf{x})=\sum_{\pi\in\mathfrak{S}} g(\pi(\mathbf{x})).
$$ 
We have the relationships $\pi^{-1}\pi(\mathbf{x})=\mathbf{x}$, $\pi^{-1}(\mu^{\otimes n})=\mu^{\otimes n}$ and $\overline{\overline{g}}(\pi(\mathbf{x}))=\overline{\overline{g}}(\mathbf{x})$. In particular, $\overline{\overline{g}}(\mathbf{x})=\overline{\overline{g}}(\mathfrak{r}(\mathbf{x}))$. For non negative Borel function $h$ on $\rr_+^n$, we compute
$$
\dcb
\mathbb{E}^{\nu_t}[gh(\mathfrak{r})]
&=&\int_{\rr_+^n} (ga_t)(\mathbf{x})h(\mathfrak{r}(\mathbf{x})) \mu^{\otimes n}(d\mathbf{x})\\

&=&\sum_{\pi\in\mathfrak{S}}\int_{\rr_+^n} (ga_t)(\mathbf{x})h(\mathfrak{r}(\mathbf{x}))\ind_{\{x_{\pi(1)}<\ldots<x_{\pi(n)}\}} \mu^{\otimes n}(d\mathbf{x})\\
&&\mbox{ because $\mu$ is atom-free}\\

&=&\sum_{\pi\in\mathfrak{S}}\int_{\rr_+^n} (ga_t)(\mathbf{x})h(x_{\pi(1)},\ldots,x_{\pi(n)})\ind_{\{x_{\pi(1)}<\ldots<x_{\pi(n)}\}} \mu^{\otimes n}(d\mathbf{x})\\

&=&\sum_{\pi\in\mathfrak{S}}\int_{\rr_+^n} (ga_t)(\pi^{-1}\pi(\mathbf{x}))h(x_{\pi(1)},\ldots,x_{\pi(n)})\ind_{\{x_{\pi(1)}<\ldots<x_{\pi(n)}\}} \mu^{\otimes n}(d\mathbf{x})\\

&=&\sum_{\pi\in\mathfrak{S}}\int_{\rr_+^n} (ga_t)(\pi^{-1}(\mathbf{x}))h(x_{1},\ldots,x_{n})\ind_{\{x_{1}<\ldots<x_{n}\}} \mu^{\otimes n}(d\mathbf{x})\\

&=&\int_{\rr_+^n} \overline{\overline{ga_t}}(\mathbf{x})h(x_{1},\ldots,x_{n})\ind_{\{x_{1}<\ldots<x_{n}\}} \mu^{\otimes n}(d\mathbf{x})\\

&=&\sum_{\pi\in\mathfrak{S}}\int_{\rr_+^n} \frac{1}{n!}\overline{\overline{ga_t}}(\pi(\mathbf{x}))h(x_{\pi(1)},\ldots,x_{\pi(n)})\ind_{\{x_{\pi(1)}<\ldots<x_{\pi(n)}\}} \mu^{\otimes n}(d\mathbf{x})\\

&=&\sum_{\pi\in\mathfrak{S}}\int_{\rr_+^n} \frac{1}{n!}\overline{\overline{ga_t}}(\mathbf{x})h(\mathfrak{r}(\mathbf{x}))\ind_{\{x_{\pi(1)}<\ldots<x_{\pi(n)}\}} \mu^{\otimes n}(d\mathbf{x})\\

&=&\int_{\rr_+^n} \frac{1}{n!}\overline{\overline{ga_t}}(\mathbf{x})h(\mathfrak{r}(\mathbf{x})) \mu^{\otimes n}(d\mathbf{x}).

\dce
$$
In particular, if $g\equiv 1$, $$
\int_{\rr_+^n} a_t(\mathbf{x})h(\mathfrak{r}(\mathbf{x})) \mu^{\otimes n}(d\mathbf{x})
=
\int_{\rr_+^n} \frac{1}{n!}\overline{\overline{a_t}}(\mathbf{x})h(\mathfrak{r}(\mathbf{x})) \mu^{\otimes n}(d\mathbf{x}).
$$
Continuing the above computation with that property, we obtain
$$
\dcb
\mathbb{E}^{\nu_t}[gh(\mathfrak{r})]
&=&
\int_{\rr_+^n}\frac{1}{n!}\overline{\overline{a_t}}(\mathbf{x}) \frac{\overline{\overline{ga_t}}(\mathfrak{r}(\mathbf{x}))}{\overline{\overline{a_t}}(\mathfrak{r}(\mathbf{x}))}\ind_{\{\overline{\overline{a_t}}(\mathfrak{r}(\mathbf{x}))>0\}}h(\mathfrak{r}(\mathbf{x})) \mu^{\otimes n}(d\mathbf{x})\\

&=&\int_{\rr_+^n}a_t(\mathbf{x}) \frac{\overline{\overline{ga_t}}(\mathfrak{r}(\mathbf{x}))}{\overline{\overline{a_t}}(\mathfrak{r}(\mathbf{x}))}\ind_{\{\overline{\overline{a_t}}(\mathfrak{r}(\mathbf{x}))>0\}}h(\mathfrak{r}(\mathbf{x})) \mu^{\otimes n}(d\mathbf{x})\\

&=&\mathbb{E}^{\nu_t}[\frac{\overline{\overline{ga_t}}(\mathfrak{r})}{\overline{\overline{a_t}}(\mathfrak{r})}\ind_{\{\overline{\overline{a_t}}(\mathfrak{r})>0\}}h(\mathfrak{r})].
\dce
$$
This computation shows that$$
\dcb
\mathbb{E}^{\nu_t}[g|\mathfrak{r}=\mathbf{x}]
=\frac{\overline{\overline{ga_t}}(\mathbf{x})}{\overline{\overline{a_t}}(\mathbf{x})}\ind_{\{\overline{\overline{a_t}}(\mathbf{x})>0\}}.	
\dce
$$
Another consequence of the above computations is that, under \cref{density}, $\overline{\tau}$ satisfies also the density hypothesis with the density function with respect to $\mu^{\otimes n}$ and to $\mathcal{F}_t$, $t\in\mathbb{R}_+$, given by $\widetilde{a}_t(\mathbf{x})=\I_\seq{x_1< x_2< \dots < x_n}\overline{\overline{a_t}}(\mathbf{x})$ (cf. \cite{EJJ}).

\

\subsubsection{Computing the $\ff^{\btau_{\varrho}}$-conditional expectation}\label{s5.1}

We use the notation introduced in
the previous subsections. Let $\varrho\in \mathcal{I}(k,n)$ and $T\subset I_k$ with $j=\#T$. Note that, for $t\in\rr_+$, the $\sigma$-algebra $\mathcal{F}^{\sigma(\btau_{\varrho(T)})}_t$ is generated by $h(\btau_{\varrho(T)})$ where $h$ runs over the family of all bounded functions on $\Omega\times\rr_+^j$, $\cap_{s>t}(\mathcal{F}_s\otimes\sigma(\rr_+^j))$ measurable.

\bl\label{fcondrhoT}
For any non negative
$\cap_{s>t}(\F_s\otimes\mathcal{B}(\rr^n))$ measurable function $g$, for any bounded $\mathbb{F}^{\sigma(\btau_{\varrho(T)})}$ stopping time $U$,
$$
\mathbb{E}_\mathbb{P}[g(\btau)|\mathcal{F}^{\sigma(\btau_{\varrho(T)})}_U] =\frac{(ga_U)_{\varrho(T)}(\btau_{\varrho(T)})	}{a_{\varrho(T),U}(\btau_{\varrho(T)})}\ind_{\{a_{\varrho(T),U}(\btau_{\varrho(T)})>0\}}.
$$
If $U$ is a bounded $\mathbb{F}^{\sigma(\btau_{\varrho(T)})}$ predictable stopping time, we also have
$$
\mathbb{E}_\mathbb{P}[g(\btau)|\mathcal{F}^{\sigma(\btau_{\varrho(T)})}_{U-}] =\frac{(ga_{U-})_{\varrho(T)}(\btau_{\varrho(T)})}{a_{\varrho(T),U-}(\btau_{\varrho(T)})}\ind_{\{a_{\varrho(T),U-}(\btau_{\varrho(T)})>0\}}.
$$
\el

\begin{proof}
By monotone convergence theorem we only need to prove the lemma for $0\leq g\leq 1$. Let us show firstly the lemma for $U=t\in\rr_+$. For a bounded function $h$ on $\Omega\times\rr_+^j$, $\cap_{s>t}(\mathcal{F}_s\otimes\sigma(\rr_+^j))$ measurable, according Lemma \ref{ahmu},$$
\begin{array}{lll}
&&\mathbb{E}_\mathbb{P}[g(\btau) h(\btau_{\varrho(T)})]\\
&=&\mathbb{E}_\mathbb{P}[\int_{\rr_+^n}g(\mathbf{x}) h(\mathbf{x}_{\varrho(T)})a_t(\mathbf{x})\mu^{\otimes n}(d\mathbf{x})]\\
&=&\mathbb{E}_\mathbb{P}[\frac{(ga_t)_{\varrho(T)}(\btau_{\varrho(T)})}{a_{\varrho(T),t}(\btau_{\varrho(T)})}\ind_{\{a_{\varrho(T),t}(\btau_{\varrho(T)})>0\}} h(\btau_{\varrho(T)})],
\end{array}
$$
where the last equality comes from Lemma \ref{ahmu} applied with respect to $\btau_{\varrho(T)}$. The formula is proved for $U=t$, because $\frac{(ga_t)_{\varrho(T)}(\btau_{\varrho(T)})}{a_{\varrho(T),t}(\btau_{\varrho(T)})}\ind_{\{a_{\varrho(T),t}(\btau_{\varrho(T)})>0\}}$ is $\mathcal{F}^{\sigma(\btau_{\varrho(T)})}_t$ measurable.

For any $n\in\mathbb{N}$, for $(\omega,\mathbf{x})\in\Omega\times\rr_+^n$, let $R^n(\omega,\mathbf{x})=\inf\{s\in\mathbb{Q}_+: a_s(\omega,\mathbf{x})>n\}$. Then, for $b\in\rr_+$, 
$$
\{R^n\geq b\}=\{(\omega,\mathbf{x})\in\Omega\times\rr_+^n: \forall s\in\mathbb{Q}_+\cap[0,b), a_s(\omega,\mathbf{x})\leq n\}
\in\mathcal{F}_{b-}\otimes\mathcal{B}(\rr_+^n).
$$
Applying the above formula at constant time $t$ to $g(\btau)\ind_{\{t< R^n(\btau)\}}$, we can write$$
\mathbb{E}_\mathbb{P}[g(\btau)\ind_{\{t< R^n(\btau)\}}|\mathcal{F}^{\sigma(\btau_{\varrho(T)})}_t] =\frac{(g\ind_{\{t< R^n\}}a_t)_{\varrho(T)}(\btau_{\varrho(T)})}{a_{\varrho(T),t}(\btau_{\varrho(T)})}\ind_{\{a_{\varrho(T),t}(\btau_{\varrho(T)})>0\}}.
$$
Note that $
g\ind_{\{t< R^n\}}a_t\leq n.
$
By the dominated convergence theorem, for almost all $\omega$, the map $$
(g\ind_{\{t< R^n\}}a_t)_{\varrho(T)}(\btau_{\varrho(T)}), t\in\rr_+,
$$ 
is right continuous. By \cite[Lemme(1.8)]{JJ}, $\ind_{\{a_{\varrho(T),t}(\btau_{\varrho(T)})>0\}}$ also is right continuous. Hence, the above formula can be extended to any bounded $\mathbb{F}^{\sigma(\btau_{\varrho(T)})}$ stopping time $U$:
$$
\mathbb{E}_\mathbb{P}[g(\btau)\ind_{\{U< R^n(\btau)\}}|\mathcal{F}^{\sigma(\btau_{\varrho(T)})}_U] =\frac{(g\ind_{\{U< R^n\}}a_U)_{\varrho(T)}(\btau_{\varrho(T)})}{a_{\varrho(T),U}(\btau_{\varrho(T)})}\ind_{\{a_{\varrho(T),U}(\btau_{\varrho(T)})>0\}}.
$$
Note that, by dominated convergence theorem, $(g\ind_{\{t< R^n\}}a_t)_{\varrho(T)}(\btau_{\varrho(T)})$ has left limit $(g\ind_{\{t\leq R^n\}}a_{t-})_{\varrho(T)}(\btau_{\varrho(T)})$. If $U$ is a bounded $\mathbb{F}^{\sigma(\btau_{\varrho(T)})}$ predictable stopping time, we also have
$$
\mathbb{E}_\mathbb{P}[g(\btau)\ind_{\{U\leq R^n(\btau)\}}|\mathcal{F}^{\sigma(\btau_{\varrho(T)})}_{U-}] =\frac{(g\ind_{\{U\leq R^n\}}a_{U-})_{\varrho(T)}(\btau_{\varrho(T)})}{a_{\varrho(T),U-}(\btau_{\varrho(T)})}\ind_{\{a_{\varrho(T),U-}(\btau_{\varrho(T)})>0\}}.
$$
Now let $n\uparrow\infty$ we prove the lemma. 
\end{proof}

Notice that in the above formulas, we can remove the indicator $\ind_{\{a_{\varrho(T),t}(\btau_{\varrho(T)})>0\}}$, because by \cite[Corollaire(1.11)]{JJ} the process $\ind_{\{a_{\varrho(T),t}(\btau_{\varrho(T)})=0\}}$ is evanescent. 

\bcor\label{essc}
For any Borel function $g$ on $\rr_+^n$ such that $g(\tau)$ is integrable, the process $(ga_t)_{\varrho(T)}(\btau_{\varrho(T)}), t\in\rr_+$, is c\`adl\`ag whose left limit is the process $(ga_{t-})_{\varrho(T)}(\btau_{\varrho(T)}), t\in\rr_+$.
\ecor

\brem
Corollary \ref{essc} implies in particular that $(a_t)_{\varrho(T)},t\in\rr_+,$ is c\`adl\`ag so that $(a_t)_{\varrho(T)}(\btau_{\varrho(T)}),t\in\rr_+,$ coincides with $a_{t,\varrho(T)}(\btau_{\varrho(T)}),t\in\rr_+$. It is an important property in practice (for example, numerical implantation) because it gives a concrete way to compute $a_{t,\varrho(T)}(\btau_{\varrho(T)})$.
\erem

We define,$$
\begin{array}{lll}
\max \mathbf{x}_{\varrho(T)}:=\max_{i\in T} x_{\varrho(i)},\
\mathbf{x}\in\rr_+^n, \\ 
\min \mathbf{x}_{\varrho(I_k\setminus T)}:=\min_{i\in I_k\setminus T} x_{\varrho(i)},\
\mathbf{x}\in\rr_+^n,\\
\mathtt{A}_{t,T,\varrho} := \seq{\mathbf{x}\in\rr_+^n: \max \mathbf{x}_{\varrho(T)}\leq
t, \ \min \mathbf{x}_{\varrho(I_k\setminus T)}> t},\ t\in\rr_+.
\end{array}
$$

\bl\label{triburest} 
For any bounded $\mathbb{F}$ stopping time $U$ we have$$
\mathcal{F}^{\btau_\varrho}_U
=\mathcal{F}_U\vee\sigma(\btau_\varrho \nmid U
).
$$
For any
subset $T$ of $I_k$, the process $\ind_{\mathtt{A}_{t,T,\varrho}}(\btau), t\in\rr_+$, is $\mathbb{F}^{\btau_\varrho}$ optional and
\bde
\mathcal{F}^{\btau_\varrho}_U\cap\{\btau\in\mathtt{A}_{U,T,\varrho}\}
=\left(\mathcal{F}_U\vee\sigma(\btau_{\varrho} \nmid U
)\right)\cap\{\btau\in\mathtt{A}_{U,T,\varrho}\}
=\left(\mathcal{F}_U\vee\sigma(\btau_{\varrho(T)})\right)\cap\{\btau\in\mathtt{A}_{U,T,\varrho}\}. 
\ede
\el

\begin{proof}
We write$$
\dcb 
\{\btau\in\mathtt{A}_{t,T,\varrho}\} 
&=\seq{\forall i\in T, \tau_{\varrho(i)}\nmid t<\infty, \ \forall
i\in I_k\setminus T, \tau_{\varrho(i)}\nmid t=\infty},
\dce 
$$
which is a $\sigma(\btau_{\varrho}\nmid t)$-measurable set. Hence, the process $\ind_{\mathtt{A}_{t,T,\varrho}}(\btau), t\in\rr_+$, is $\mathbb{F}^{\btau_\varrho}$ adapted. But it also is c\`adl\`ag. The first assertion is proved.

Let
$\boldsymbol\gamma=(\gamma_{1},\dots,\gamma_{k})$ be the increasing re-ordering of $\btau_\varrho$ and set $\gamma_{0}=0, \gamma_{k+1}=\infty$. The density hypothesis with respect to $\mathbb{F}$ holds for
$\btau_\varrho$, since it holds for $\btau$. It
is proved in \cite{S} that the optional splitting formula holds
with respect to $\mathbb{F}^{\boldsymbol{\btau_{\varrho}}}$. As a consequence, for
$0\leq j\leq k$,
\begin{align*}
\mathcal{F}^{\btau_\varrho}_U\cap \{\gamma_{j}\leq U<
\gamma_{j+1}\}
&=\Big(\mathcal{F}_U\vee\sigma(\btau_\varrho\nmid\gamma_{j} )\Big)\cap\{\gamma_{j}\leq U< \gamma_{j+1}\}\\
&=\Big(\mathcal{F}_U\vee\sigma(\btau_\varrho \nmid U
)\Big)\cap\{\gamma_{j}\leq U< \gamma_{j+1}\}.
\end{align*}
Notice that$$
\ind_{\{\gamma_j\leq U<\gamma_{j+1}\}}=\sum_{T\subset I_k:\#T=j}\I_{\mathtt{A}_{U,T,\varrho}}(\btau)
$$
is $\mathcal{F}^{\btau_\varrho}_U$ as well as $\mathcal{F}_U\vee\sigma(\btau_\varrho \nmid U
)$ measurable. Hence, we can take the union of the above identities to conclude$$
\mathcal{F}^{\btau_\varrho}_U
=\mathcal{F}_U\vee\sigma(\btau_\varrho \nmid U
).
$$
If   the cardinal of $T$  is
equal to $j$, we have $\mathtt{A}_{U,T,\varrho}\subset\{\gamma_{j}\leq
U< \gamma_{j+1}\}$ and the last claim of the lemma follows from
the above identities together with the fact that $\sigma(\btau_{\varrho} \nmid U
)\cap\mathtt{A}_{U,T,\varrho} =\sigma(\btau_{\varrho(T)})\cap\mathtt{A}_{U,T,\varrho}$.
\end{proof}

We introduce another notations $U_{T,\varrho}(\mathbf{x}):=\max \mathbf{x}_{\varrho(T)}$ and $S_{T,\varrho}(\mathbf{x}):=\min \mathbf{x}_{\varrho(I_k\setminus T)}$. Recall (cf. \cite{J1}) that, for any random time $U$, $\F_U$ (resp. $\F_{U-}$) denotes the $\sigma$-algebra generated by $K_U$, where $K$ denotes a $\mathbb{F}$ optional (resp. predictable) process.

\bl\label{attr} 
For bounded $\mathbb{F}^{\btau_\varrho}$ stopping time $U$, for any
non negative $\F_U\otimes\mathcal{B}(\rr^n)$ measurable function $g$, we have 
\bde
\mathbb{E}_\mathbb{P}(g(\btau)|\mathcal{F}^{\btau_\varrho}_U) =
\sum_{T\subset I_k}\I_{\seq{U_{T,\varrho}(\btau)\leq
U< S_{T,\varrho}(\btau)}}
\frac{(\ind_{\{U<S_{T,\varrho}\}}ga_U)_{\varrho(T)}(\btau_{\varrho(T)})}{(\ind_{\{U<S_{T,\varrho}\}}a_U)_{\varrho(T)}(\btau_{\varrho(T)})}. 
\ede
\el 

\begin{proof}
By monotone convergence theorem we only need to prove the lemma for bounded $g$. By monotone class theorem, we only need to prove the lemma for bounded Borel function $g$ on $\rr_+^n$. Let us firstly consider $U=t\in\rr_+$. The parameters $U$, $T$ and $\varrho$ being fixed, for simplicity we write $\mathtt{A}$ instead of $\mathtt{A}_{t,T,\varrho}$. We also introduce
$$
\mathtt{F}=\{\mathbf{x}\in\rr_+^n:\min \mathbf{x}_{\varrho(I_k\setminus T)}> t\}.
$$ 
For any bounded
$\mathcal{B}(\rr^{k})$-measurable function $h$, there exists a $\mathcal{B}(\rr^{\# T})$-measurable function $h'$ such that $h(\btau_{\varrho}\nmid t)=h'(\btau_{\varrho(T)})$ on $\mathtt{A}$. Let $B\in\mathcal{F}_t$. With help of Lemma \ref{fcondrhoT}, we compute
$$
\begin{array}{lll}
&\mathbb{E}_\mathbb{P}\left(\I_B h(\btau_{\varrho}\nmid t)\I_{\mathtt{A}}(\btau)\mathbb{E}[g(\btau)|\mathcal{F}^{\btau_{\varrho}}_t]\right)\\

=&\mathbb{E}_\mathbb{P}\left(\I_B h'(\btau_{\varrho(T)})\I_{\mathtt{A}}(\btau)g(\btau)\right)\\

=&\mathbb{E}_\mathbb{P}\left(\I_B h'(\btau_{\varrho(T)})\I_{\{\max \btau_{\varrho(T)}\leq t\}}\I_{\{\min \btau_{\varrho(I_k\setminus T)}> t\}}g(\btau)\right)\\

=&\mathbb{E}_\mathbb{P}\left(\I_B h'(\btau_{\varrho(T)})\I_{\{\max \btau_{\varrho(T)}\leq t\}}\mathbb{E}_\mathbb{P}[\I_{\{\min \btau_{\varrho(I_k\setminus T)}> t\}}g(\btau)|\mathcal{F}^{\sigma(\btau_{\varrho(T)})}_t]\right)\\

=&\mathbb{E}_\mathbb{P}\left(\I_B h'(\btau_{\varrho(T)})\I_{\{\max \btau_{\varrho(T)}\leq t\}}\frac{(\ind_{\mathtt{F}}ga_t)_{\varrho(T)}(\btau_{\varrho(T)})}{(a_t)_{\varrho(T)}(\btau_{\varrho(T)})}\ind_{\{(a_t)_{\varrho(T)}(\btau_{\varrho(T)})>0\}}\ind_{\{(\ind_{\mathtt{F}}a_t)_{\varrho(T)}(\btau_{\varrho(T)})>0\}}\right)\\

=&\mathbb{E}_\mathbb{P}\left(\I_B h'(\btau_{\varrho(T)})\I_{\{\max \btau_{\varrho(T)}\leq t\}}\I_{\{\min \btau_{\varrho(I_k\setminus T)}> t\}}\right.\\

&\left.\frac{(a_t)_{\varrho(T)}(\btau_{\varrho(T)})}{(\ind_{\mathtt{F}}a_t)_{\varrho(T)}(\btau_{\varrho(T)})}\frac{(\ind_{\mathtt{F}}ga_t)_{\varrho(T)}(\btau_{\varrho(T)})}{(a_t)_{\varrho(T)}(\btau_{\varrho(T)})}\ind_{\{(a_t)_{\varrho(T)}(\btau_{\varrho(T)})>0\}}\ind_{\{(\ind_{\mathtt{F}}a_t)_{\varrho(T)}(\btau_{\varrho(T)})>0\}}\right)\\

=&\mathbb{E}_\mathbb{P}\left(\I_B h(\btau_{\varrho}\nmid t)\I_{\mathtt{A}}(\btau)
\frac{(\ind_{\mathtt{F}}ga_t)_{\varrho(T)}(\btau_{\varrho(T)})}{(\ind_{\mathtt{F}}a_t)_{\varrho(T)}(\btau_{\varrho(T)})}\right).\\

\end{array}
$$
It is to note that the random variables$$
\I_{\mathtt{A}}(\btau)
\mbox{ and }
\I_{\mathtt{A}}(\btau)
\frac{(\ind_{\mathtt{F}}ga_t)_{\varrho(T)}(\btau_{\varrho(T)})}{(\ind_{\mathtt{F}}a_t)_{\varrho(T)}(\btau_{\varrho(T)})}
$$
are $\mathcal{F}^{\btau_\varrho}_t=\mathcal{F}_t\vee\sigma(\btau_{\varrho}\nmid t)$ measurable. By Lemma \ref{triburest}), the above computation implies that$$
\I_{\mathtt{A}}(\btau)\mathbb{E}[g(\btau)|\mathcal{F}^{\btau_{\varrho}}_t]
=
\I_{\mathtt{A}}(\btau)
\frac{(\ind_{\mathtt{F}}ga_t)_{\varrho(T)}(\btau_{\varrho(T)})}{(\ind_{\mathtt{F}}a_t)_{\varrho(T)}(\btau_{\varrho(T)})}.
$$
Recall that $\boldsymbol\gamma=(\gamma_{1},\dots,\gamma_{k})$ is the increasing re-ordering of $\btau_\varrho$ and $\gamma_{0}=0, \gamma_{k+1}=\infty$. Notice also that, under \cref{density}, 
$\P(\tau_i = \tau_j)=0$ for any pair of $i,j$ such that $i \neq
j$. It results that$$
\sum_{T\subset I_k:\#T=j}\I_{\mathtt{A}_{t,T,\varrho}}(\btau)
=\ind_{\{\gamma_j\leq t<\gamma_{j+1}\}},
$$
and 
$$
\dcb
\mathbb{E}[g(\btau)|\mathcal{F}^{\btau_{\varrho}}_t]
&=&\sum_{j=0}^k\ind_{\{\gamma_j\leq t<\gamma_{j+1}\}}\mathbb{E}[g(\btau)|\mathcal{F}^{\btau_{\varrho}}_t]\\
&=&\sum_{j=0}^k\sum_{T\subset I_k:\#T=j}\I_{\mathtt{A}_{t,T,\varrho}}(\btau)\mathbb{E}[g(\btau)|\mathcal{F}^{\btau_{\varrho}}_t]\\
&=&\sum_{T\subset I_k}
\I_{\seq{U_{T,\varrho}(\btau)\leq
t< S_{T,\varrho}(\btau)}}
\frac{(\ind_{\{t<S_{T,\varrho}\}}ga_t)_{\varrho(T)}(\btau_{\varrho(T)})}{(\ind_{\{t<S_{T,\varrho}\}}a_t)_{\varrho(T)}(\btau_{\varrho(T)})}.
\dce
$$
Applying Corollary \ref{essc}, we extend this formula to any bounded $\mathbb{F}^{\btau_{\varrho}}$ stopping time $U$.
\end{proof}

\bcor\label{coattr} 
For any bounded $\mathbb{F}^{\btau_{\varrho}}$ predictable stopping time $U$, for any
non negative $\F_{U-}\otimes \mathcal{B}(\rr^n)$ measurable function $g$, we have 
\bde
\mathbb{E}_\mathbb{P}(g(\btau)|\mathcal{F}^{\btau_\varrho}_{U-}) =
\sum_{T\subset I_k}
\I_{\seq{U_{T,\varrho}(\btau)<
U\leq S_{T,\varrho}(\btau)}}
\frac{(\ind_{\{U\leq S_{T,\varrho}\}}ga_{U-})_{\varrho(T)}(\btau_{\varrho(T)})}{(\ind_{\{U\leq S_{T,\varrho}\}}a_{U-})_{\varrho(T)}(\btau_{\varrho(T)})}. 
\ede
\ecor

\

\subsubsection{Computing the $\ff^{\overline{\btau}_{k}}$-conditional expectation}\label{s5.2}

The following lemma is straightforward.

\bl For $t\in\rr_+$, for any non negative
$\F_t\otimes\mathcal{B}(\rr^n)$ measurable function $g$, we have 
\bde
\mathbb{E}_\mathbb{P}[g(\btau)|\mathcal{F}^{\sigma(\overline{\btau})}_t] 
=\mathbb{E}^{\nu_t}[g|\mathfrak{r}=\overline{\btau}]
=\frac{\overline{\overline{ga_t}}(\overline{\btau})}{\overline{\overline{a_t}}(\overline{\btau})}\ind_{\{\overline{\overline{a_t}}(\overline{\btau})>0\}}.
\ede 
\el

We will denote the last random variable by $\ddot{g}(\overline{\btau})$. We now can apply the results in subsection \ref{s5.1} on the vector $\overline{\btau}$ with $\varrho$ being the identity map $\iota_k$ in $I_k$. Then, $
\ind_{\mathtt{A}_{t,T,\iota_k}}(\overline{\btau})\neq 0,
$
only if $T$ is of the form $T=I_j=\{1,\ldots,j\}$ for $0\leq j\leq k$ ($T=\emptyset$ if $j=0$) and, in this case,$$
\ind_{\mathtt{A}_{t,I_j,\iota_k}}(\overline{\btau})
=
\ind_{\{\tau_{(j)}\leq t<\tau_{(j+1)}\}},\ 
\ind_{\mathtt{F}_{t,I_j,\iota_k}}(\overline{\btau})
=
\ind_{\{t<\tau_{(j+1)}\}},\ 
\mbox{ and }
\overline{\btau}_{\iota_k(T)}=\overline{\btau}_j.
$$
Let $\mathfrak{p}_i$ be the projection $\mathfrak{p}_i(\mathbf{x})=x_i$. From Lemma \ref{attr} we obtain

\bl\label{attr2} 
For any bounded $\mathbb{F}^{\overline{\btau}_k}$ stopping time $U$, for any non negative
$\F_U\otimes\mathcal{B}(\rr^n)$ measurable function $g$, we have 
\bde
\mathbb{E}_\mathbb{P}(g(\btau)|\mathcal{F}^{\overline{\btau}_k}_U) =
\sum_{j=0}^k\ind_{\{\tau_{(j)}\leq U<\tau_{(j+1)}\nmid\tau_{(k)}\}}
\frac{(\ind_{\{U<\mathfrak{p}_{j+1}\}}\ddot{g}\widetilde{a}_U)_{I_j}(\overline{\btau}_j)}{(\ind_{\{U<\mathfrak{p}_{j+1}\}}\widetilde{a}_U)_{I_j}(\overline{\btau}_{j})}. 
\ede
For any bounded $\mathbb{F}^{\overline{\btau}_k}$ predictable stopping time $U$, for any non negative
$\F_{U-}\otimes\mathcal{B}(\rr^n)$ measurable function $g$, we have 
\bde
\mathbb{E}_\mathbb{P}(g(\btau)|\mathcal{F}^{\overline{\btau}_k}_{U-}) =
\sum_{j=0}^k\ind_{\{\tau_{(j)}< t\leq \tau_{(j+1)}\nmid\tau_{(k)}\}}
\frac{(\ind_{\{U\leq \mathfrak{p}_{j+1}\}}\ddot{g}\widetilde{a}_{U-})_{I_j}(\overline{\btau}_j)}{(\ind_{\{U\leq\mathfrak{p}_{j+1}\}}\widetilde{a}_{U-})_{I_j}(\overline{\btau}_{j})}. 
\ede
\el

\

\subsubsection{The $\ff^{\btau_{\varrho}}$ drift computation}\label{btauvarrho}

We fix a $\varrho\in \mathcal{I}(k,n)$. We recall the notations $T\subset I_k$,  $\mathtt{A}_{s,T,\varrho}
:=\seq{\max \btau_{\varrho(T)}\leq s, \min \btau_{\varrho(I_k\setminus T)} > s}$, $U_{T,\varrho}(\mathbf{x}):=\max \mathbf{x}_{\varrho(T)}$ and $S_{T,\varrho}(\mathbf{x}):=\min \mathbf{x}_{\varrho(I_k\setminus T)}$, also $(\gamma_1,\dots,
\gamma_k)$ representing the increasing re-ordering of $(\tau_{(\varrho(1))},\dots,
\tau_{(\varrho(k))})$.

In next subsection we will compute the $\ff^{\overline{\btau}_{k}}$-semimartingale decomposition of a bounded $\ff$-martingale $M$ by first computing $\I_{D_\varrho}\wh V^\varrho$, which is the stochastic integral of $\frac{\I_{D_\varrho}}{N^\varrho_-}$ against the $\ff^{\btau_\varrho}$-drift of $N^\varrho M$ and then computing its dual $\ff^{\overline{\btau}_{k}}$-predictable projection. This subsection is devoted to the computation of the drift of $N^\varrho M$. We notice that this computation is an extension of the usual $\ff^{\btau_\varrho}$ semimartingale decomposition formula for $M$.

For a bounded Borel function $g$ on $\rr^n_+$ we set $$
L^g_t =
\cE{\P}{g(\btau)}{\F^{\btau_\varrho}_t},\ t\in\rr_+.
$$ 
A formula is given in Lemma \ref{attr} dealing with $L^g$. But that formula is not adapted to the computation that we will do in this subsection. From \cite[Th\'eor\`eme (2.5)]{JJ}, there exists a $\mathcal{P}(\ff)\otimes
\mathcal{B}(\rr^n_+)$-measurable process $u^M_v(\omega, \mathbf{x})$ such
that 
\begin{equation}\label{uvjacod}
d\pbracket{a(\mathbf{x})}{M}_v = u^M_v(\mathbf{x})d\pbracket{M}{M}_v,
\end{equation} 
where the predictable bracket $\pbracket{M}{M}$ is calculated
in the filtration $\ff$. The process $u(\mathbf{x})$ is known to satisfy
\begin{equation}\label{ufini}
\int_0^t \frac{1}{a_{s-}(\btau)}|u^M_s(\btau)|d\cro{M,M}_s<\infty, \ \forall t\in\rr_+,
\end{equation}
(so that we assume that $u^M_s(\mathbf{x})=0$ whenever $a_{s-}(\mathbf{x})=0$). But the computations below will require a stronger condition. 

\bhyp\label{aAA}
There exists an increasing sequence $(R_n)_{n\in\mathbb{N}}$ of bounded $\mathbb{F}$ stopping times such that $\sup_{n\in\nn_+}R_n=\infty$ and 
$$
\mathbb{E}[\int_0^{R_n} \frac{1}{a_{s-}(\btau)}|u^M_s(\btau)|d\cro{M,M}_s]<\infty, \ \forall n\in\nn_+.
$$
\ehyp
Notice that the above inequality is equivalent to 
$$
\mathbb{E}[\int \int_0^{R_n} |u^M_s(\mathbf{x})|d\cro{M,M}_s\mu^{\otimes n}(d\mathbf{x})]<\infty.
$$
We give here a sufficient condition for Assumption \ref{aAA} to hold for any bounded $\mathbb{F}$ martingale $M$.

\bl\label{umum}
Suppose that there exists an increasing sequence $(R_n)_{n\in\mathbb{N}}$ of bounded $\mathbb{F}$ stopping times such that $$
\int \mathbb{E}[\sqrt{[a(\mathbf{x}),a(\mathbf{x})]_{R_n}}]\mu^{\otimes n}(d\mathbf{x})<\infty.
$$ 
Then, for any bounded $\mathbb{F}$ martingale $M$,$$
\mathbb{E}[\int \int_0^{R_n} |u^M_s(\mathbf{x})|d\cro{M,M}_s\mu^{\otimes n}(d\mathbf{x})]<\infty
$$
for all $n\in\mathbb{N}$, i.e. Assumption \ref{aAA} holds.
\el

\begin{proof}
Let $H_s(\mathbf{x})=\mathtt{sign}u_s(\mathbf{x})$. We have
$$
\dcb
&&\int \mathbb{E}[\int_0^{R_n}\left|d[a(\mathbf{x}),H(\mathbf{x})\centerdot M]_s\right|]\mu^{\otimes n}(d\mathbf{x})\\
&=&\int \sqrt{2}\mathbb{E}[[a(\mathbf{x}),a(\mathbf{x})]^{1/2}_{R_n}]\|H(\mathbf{x})\centerdot M\|_{\mathtt{BMO}}\mu^{\otimes n}(d\mathbf{x}).

\dce
$$
We note that $\|H(\mathbf{x})\centerdot M\|_{\mathtt{BMO}}$ is computed by its bracket (cf. \cite[Theorem 10.9]{HWY}) so that it is uniformly bounded by a multiple of $\|M_\infty\|_\infty$. This boundedness together with the assumption of the lemma enables us to write
$$
\dcb
&&\mathbb{E}[\int \int_0^{R_n} |u^M_s(\mathbf{x})|d\cro{M,M}_s\mu^{\otimes n}(d\mathbf{x})]\\

&=&\int \mathbb{E}[\int_0^{R_n} H_s(\mathbf{x})u^M_s(\mathbf{x}) d\cro{M,M}_s]\mu^{\otimes n}(d\mathbf{x})\\

&=&\int \mathbb{E}[\int_0^{R_n} d\cro{a(\mathbf{x}),H(\mathbf{x})\centerdot M}_s]\mu^{\otimes n}(d\mathbf{x})\\

&=&\int \mathbb{E}[[a(\mathbf{x}),H(\mathbf{x})\centerdot M]_{R_n}]\mu^{\otimes n}(d\mathbf{x})<\infty.
\dce
$$
\end{proof}

\bt\label{drifttauvarrho}
Under Assumption \ref{aAA}, for any bounded Borel function $g$ on $\rr^n_+$, the drift of $L^gM$ in $\mathbb{F}^{\btau_{\varrho}}$ is given by$$
\dcb
&&
\sum_{T\subset I_k}\int_0^{t}\ind_{\{U_{T,\varrho}(\btau_\varrho)< v\leq S_{T,\varrho}(\btau_\varrho)\}}\frac{(\ind_{\{v\leq S_{T,\varrho}\}}gu^M_v)_{\varrho(T)}(\btau_{\varrho(T)})}{(\ind_{\{v\leq S_{T,\varrho}\}}a_{v-})_{\varrho(T)}(\btau_{\varrho(T)})}d\cro{M,M}_v,\ t\in\rr_+.
\dce
$$
\et

\begin{proof}
Let $T\subset I_k$ and $j=\#T$. Let $R$ be one of $R_n$ in Assumption \ref{aAA}. We compute the following for $s,t\in\rr_+, s\leq t,$ $B \in \F_s$ and $h$ a bounded Borel function on $\rr_+^j$. Note that Fubini's theorem can be applied because of, on the one hand, the boundedness of $M$ and, on the other hand, of Assumption \ref{aAA}. As $\varrho$, $T$, and $t\in\rr_+$ are known, we simply write $U(\mathbf{x})=U_{T,\varrho}(\mathbf{x})$ and $S(\mathbf{x})=S_{T,\varrho}(\mathbf{x})$. 
$$
\dcb
&&\mathbb{E}_\mathbb{P}(\int_0^\infty\I_Bh(\btau_{\varrho(T)})\ind_{\{s<v\leq t\}}\ind_{\{U(\btau_{\varrho})<v\leq S(\btau_{\varrho})\}}\ind_{\{v\leq R\}}d(L^gM)_v)\\

&=&\mathbb{E}_\mathbb{P}(\I_Bh(\btau_{\varrho(T)})\ind_{\{s\vee U(\btau_{\varrho})<t\wedge S(\btau_{\varrho})\wedge R\}}\left((L^gM)_{t\wedge S(\btau_{\varrho})\wedge R} - (L^gM)_{(s\vee U(\btau_{\varrho}))\wedge (t\wedge S(\btau_{\varrho})\wedge R)}\right))\\

&=&\mathbb{E}_\mathbb{P}(\I_Bh(\btau_{\varrho(T)})\ind_{\{s\vee U(\btau_{\varrho})<t\wedge S(\btau_{\varrho})\wedge R\}}g(\btau)\left(M_{t\wedge S(\btau_{\varrho})\wedge R} - M_{(s\vee U(\btau_{\varrho}))\wedge (t\wedge S(\btau_{\varrho})\wedge R)}\right))\\

&=&\mathbb{E}_\mathbb{P}(\I_Bh(\btau_{\varrho(T)})g(\btau)\left(M_{t\wedge S(\btau_{\varrho})\wedge R} - M_{(s\vee U(\btau_{\varrho}))\wedge (t\wedge S(\btau_{\varrho})\wedge R)}\right))\\

&=&\mathbb{E}_\mathbb{P}[ \ind_B\int h(\mathbf{x}_{\varrho(T)})g(\mathbf{x})(M_{t\wedge S(\mathbf{x}_\varrho)\wedge R}-M_{(s\vee U(\mathbf{x}_\varrho))\wedge (t\wedge S(\mathbf{x}_\varrho)\wedge R)})a_t(\mathbf{x})\mu^{\otimes n}(d\mathbf{x})]\\

&=&\int h(\mathbf{x}_{\varrho(T)})g(\mathbf{x})\mathbb{E}_\mathbb{P}[ \ind_B(M_{t\wedge S(\mathbf{x}_\varrho)\wedge R}-M_{(s\vee U(\mathbf{x}_\varrho))\wedge (t\wedge S(\mathbf{x}_\varrho)\wedge R)})a_t(\mathbf{x})]\mu^{\otimes n}(d\mathbf{x})\\

&=&\int h(\mathbf{x}_{\varrho(T)})g(\mathbf{x})\mathbb{E}_\mathbb{P}[ \ind_B\int_{(s\vee U(\mathbf{x}_\varrho))}^{t\wedge S(\mathbf{x}_\varrho)\wedge R} u^M_v(\mathbf{x})d\cro{M,M}_v]\mu^{\otimes n}(d\mathbf{x})\\
	
&=&\int h(\mathbf{x}_{\varrho(T)})g(\mathbf{x})\mathbb{E}_\mathbb{P}[ \ind_B\int_s^t
\ind_{\{U(\mathbf{x}_\varrho)<v\leq S(\mathbf{x}_\varrho)\wedge R\}} u^M_v(\mathbf{x})d\cro{M,M}_v]\mu^{\otimes n}(d\mathbf{x})\\

&=&\mathbb{E}_\mathbb{P}[\ind_B  \int_s^{t}\left(\int h(\mathbf{x}_{\varrho(T)}) \ind_{\{U(\mathbf{x}_\varrho)<v\leq S(\mathbf{x}_\varrho)\wedge R\}} g(\mathbf{x})u^M_v(\mathbf{x})\mu^{\otimes n}(d\mathbf{x})\right)d\cro{M,M}_v]\\

&=&\mathbb{E}_\mathbb{P}[\ind_B  \int_s^{t}\left(\int h(\mathbf{x}_{\varrho(T)}) \ind_{\{U(\mathbf{x}_\varrho)<v\leq S(\mathbf{x}_\varrho)\wedge R\}} \frac{(\ind_{\{v\leq S\}}gu^M_v)_{\varrho(T)}(\mathbf{x}_{\varrho(T)})}{(\ind_{\{v\leq S\}}a_{v-})_{\varrho(T)}(\mathbf{x}_{\varrho(T)})}a_{v-}(\mathbf{x})\mu^{\otimes n}(d\mathbf{x})\right)d\cro{M,M}_v]\\

&=&\mathbb{E}_\mathbb{P}(\I_B \int_s^{t} h(\btau_{\varrho(T)})\ind_{\{U(\btau_\varrho)<v\leq S(\btau_\varrho)\wedge R\}}\frac{(\ind_{\{v\leq S\}}gu^M_v)_{\varrho(T)}(\btau_{\varrho(T)})}{(\ind_{\{v\leq S\}}a_{v-})_{\varrho(T)}(\btau_{\varrho(T)})}d\cro{M,M}_v)\\
&&\mbox{ consequence of \cite[Lemme(1.10)]{JJ},}\\

&=&\mathbb{E}_\mathbb{P}(\int_0^{\infty}\I_Bh(\btau_{\varrho(T)})\ind_{\{s<v\leq t\}}\ind_{\{U(\btau_\varrho)<v\leq S(\btau_\varrho)\}}\ind_{\{v\leq R\}}
\frac{(\ind_{\{v\leq S\}}gu^M_v)_{\varrho(T)}(\btau_{\varrho(T)})}{(\ind_{\{v\leq S\}}a_{v-})_{\varrho(T)}(\btau_{\varrho(T)})}d\cro{M,M}_v).
\dce
$$
Note that $R=R_n$ tends to the infinity and the processes $\I_Bh(\btau_{\varrho(T)})\ind_{\{s<v\leq t\}}$ generate all bounded $\mathbb{F}^{\btau_{\varrho}}$ predictable processes on $(U(\btau_	{\varrho}),S(\btau_{\varrho})]$. The above computation means that the drift of $\ind_{(U(\btau_{\varrho}),S(\btau_{\varrho})]}\centerdot(L^gM)$ is given by
$$
\int_0^{t}\ind_{\{U(\btau_\varrho)< v\leq S(\btau_\varrho)\}}\frac{(\ind_{\{v\leq S\}}gu^M_v)_{\varrho(T)}(\btau_{\varrho(T)})}{(\ind_{\{v\leq S\}}a_{v-})_{\varrho(T)}(\btau_{\varrho(T)})}d\cro{M,M}_v,\ t\in\rr_+.
$$
The lemma follows because $$
\sum_{T\subset I_k}\ind_{(U_{T,\varrho}(\btau_\varrho),S_{T,\varrho}(\btau_\varrho)]}
=\ind_{(0,\infty]}.
$$
\end{proof}

\

\subsubsection{The $\ff^{\overline{\btau}_{k}}$-semimartingale decomposition}\label{gfinal}

Denote by $\mathfrak{p}(\varrho)_{i}$ the map $\mathbf{x}\rightarrow\mathfrak{p}_{i}(\mathbf{x}_{\varrho})$.

\bt\label[theorem]{finalT}
Under Assumption \ref{aAA}, the $\ff^{\overline{\btau}_{k}}$-drift of the bounded $\ff$-martingale $M$ is given by
\bde
\sum_{\varrho \in \mathcal{I}(k,n)}\sum_{j=0}^k\int_0^{t}\ind_{\{\tau_{(j)}< v\leq\tau_{(j+1)}\nmid\tau_{(k)}\}}
\frac{(\ind_{\{v\leq \mathfrak{p}(\varrho)_{j+1}\}}\zeta_\varrho u^M_v)_{\varrho(I_j)}(\overline{\btau}_{j}))}{(\ind_{\{v\leq \mathfrak{p}(\varrho)_{j+1}\}}\zeta_\varrho a_{v-})_{\varrho(I_j)}(\overline{\btau}_{j}))}
\frac{(\ind_{\{v\leq \mathfrak{p}_{j+1}\}}\widetilde{a}^\varrho_{v-})_{I_j}(\overline{\btau}_{j})}{(\ind_{\{v\leq \mathfrak{p}_{j+1}\}}\widetilde{a}_{v-})_{I_j}(\overline{\btau}_{j})}
d\cro{M,M}_v,
\ede
for $t\in\rr_+$, where $\zeta_{\varrho}(\mathbf{x})=\ind_{\{x_{\varrho(1)}<\ldots<x_{\varrho(k)}<\min\mathbf{x}_{I_n\setminus\varrho(I_k)}\}}$, and $$
\widetilde{a}^\varrho_{v-}(\mathbf{x})
=
\ind_{\{x_1<x_2<\ldots<x_n\}}\sum_{\pi\in\mathfrak{S}: \pi(\varrho(i))=i, \forall i\in I_k)}a_{v-}(\pi(\mathbf{x})),\ \mathbf{x}\in\rr_+^n.
$$
\et

\begin{proof}
According to Lemma \ref{multiprop}, we only need to calculate the
$\ff^{\overline{\btau}_{k}}$-predictable process of finite
variation $\wt N^\varrho_-\ast\psi(\wh V^\varrho)$ for $\varrho \in \mathcal{I}(k,n)$. Note that $\zeta_\varrho(\btau)=\I_{d_\varrho}=\I_{D_\varrho}$ because of \cref{density}. On the set $\{x_{\varrho(1)}<\ldots<x_{\varrho(k)}<\min\mathbf{x}_{I_n\setminus\varrho(I_k)}\}$, $\{\mathbf{x}: U_{T,\varrho}(\mathbf{x}_\varrho)< v\leq S_{T,\varrho}(\mathbf{x}_\varrho)\}=\emptyset$ if $T\neq I_j$ with $j=\#T$, while $\{\mathbf{x}:U_{I_j,\varrho}(\mathbf{x}_\varrho)< v\leq S_{I_j,\varrho}(\mathbf{x}_\varrho)\}=\{\mathbf{x}:x_{\varrho(j)}<v\leq x_{\varrho(j+1)}\}$.

If we set $g=\zeta_\varrho$ in the above Theorem \ref{drifttauvarrho}, we obtain the drift of $N^{\varrho}M$ in $\mathbb{F}^{\btau_{\varrho}}$. According to Lemma \ref{gdel}, this drift process coincides with $N^\varrho_-\centerdot \widehat{V}^\varrho$. Consequently,$$
\ind_{D_\varrho}\widehat{V}^\varrho_t
=
\sum_{j=0}^k\int_0^{t}\frac{\zeta_\varrho(\btau)}{N^\varrho_{v-}}\ind_{\{\tau_{\varrho(j)}<v\leq \tau_{\varrho(j+1)}\}}\frac{(\ind_{\{v\leq \mathfrak{p}(\varrho)_{j+1}\}}\zeta_\varrho u^M_v)_{\varrho(I_j)}(\btau_{\varrho(I_j)})}{(\ind_{\{v\leq S_{I_j,\varrho}\}}a_{v-})_{\varrho(I_j)}(\btau_{\varrho(I_j)})}d\cro{M,M}_v,
$$
for $t\in\rr_+$. Note that $N^\varrho_{v-}$ is computed by Corollary \ref{coattr}. On the set $d_\varrho$,
$$
\dcb
N^\varrho_{v-}
&=&\mathbb{E}_\mathbb{P}(\zeta_\varrho(\btau)|\mathcal{F}^{\btau_\varrho}_{v-})\\

&=&
\sum_{T\subset I_k}\I_{\seq{U_{T,\varrho}(\btau)<
v\leq S_{T,\varrho}(\btau)}}
\frac{(\ind_{\{v\leq S_{T,\varrho}\}}\zeta_\varrho a_{v-})_{\varrho(T)}(\btau_{\varrho(T)})}{(\ind_{\{v\leq S_{T,\varrho}\}}a_{v-})_{\varrho(T)}(\btau_{\varrho(T)})}\\

&=&
\sum_{j=0}^k\I_{\seq{U_{I_j,\varrho}(\btau)<
v\leq S_{I_j,\varrho}(\btau)}}
\frac{(\ind_{\{v\leq \mathfrak{p}(\varrho)_{j+1}\}}\zeta_\varrho a_{v-})_{\varrho(I_j)}(\btau_{\varrho(I_j)})}{(\ind_{\{v\leq S_{I_j,\varrho}\}}a_{v-})_{\varrho(I_j)}(\btau_{\varrho(I_j)})}. 
\dce
$$
This yields
$$
\dcb
\ind_{D_\varrho}\widehat{V}^\varrho_t
&=&
\sum_{j=0}^k\int_0^{t}\zeta_\varrho(\btau)\ind_{\{\tau_{\varrho(j)}<v\leq \tau_{\varrho(j+1)}\}}\frac{(\ind_{\{v\leq \mathfrak{p}(\varrho)_{j+1}\}}\zeta_\varrho u^M_v)_{\varrho(I_j)}(\btau_{\varrho(I_j)})}{(\ind_{\{v\leq \mathfrak{p}(\varrho)_{j+1}\}}\zeta_\varrho a_{v-})_{\varrho(I_j)}(\btau_{\varrho(I_j)})}d\cro{M,M}_v\\

&=:&
\sum_{j=0}^k\int_0^{t}\ind_{\{\tau_{(j)}<v\leq \tau_{(j+1)}\}}\varphi_{v,j,\varrho}(\btau)d\cro{M,M}_v.
\dce
$$ 
By Lemma \ref{gdel} we need to compute the $\mathbb{F}^{\overline{\btau}_k}$ dual predictable projection of the above process to obtain $\wt N^\varrho_-\ast\psi(\wh V^\varrho)$. Since $\cro{M,M}$ is $\mathbb{F}$ predictable, it is enough to compute the predictable projection of the integrand. We compute firstly
$$
\dcb
&&\ddot{\varphi}_{v,j,\varrho}\widetilde{a}_{v-}(\mathbf{x})\\
&=&
\frac{\overline{\overline{\varphi_{v,j,\varrho}a_{v-}}}(\mathbf{x})}{\overline{\overline{a_{v-}}}(\mathbf{x})}\ind_{\{\overline{\overline{a_{v-}}}(\mathbf{x})>0\}}\ind_{\{x_1<x_2<\ldots<x_n\}}\overline{\overline{a_{v-}}}(\mathbf{x})\\

&=&
\overline{\overline{\varphi_{v,j,\varrho}a_{v-}}}(\mathbf{x})\ind_{\{x_1<x_2<\ldots<x_n\}}\\

&=&
\sum_{\pi\in\mathfrak{S}}
\varphi_{v,j,\varrho}(\pi(\mathbf{x}))a_{v-}(\pi(\mathbf{x}))\ind_{\{x_1<x_2<\ldots<x_n\}}\\

&=&
\sum_{\pi\in\mathfrak{S}}
\ind_{\{x_{\pi(\varrho(1))}<\ldots<x_{\pi(\varrho(k))}<\min\pi(\mathbf{x})_{I_n\setminus\varrho(I_k)}\}}\frac{(\ind_{\{v\leq \mathfrak{p}(\varrho)_{j+1}\}}\zeta_\varrho u^M_v)_{\varrho(I_j)}(\pi(\mathbf{x})_{\varrho(I_j)})}{(\ind_{\{v\leq \mathfrak{p}(\varrho)_{j+1}\}}\zeta_\varrho a_{v-})_{\varrho(I_j)}(\pi(\mathbf{x})_{\varrho(I_j)})}a_{v-}(\pi(\mathbf{x}))\ind_{\{x_1<x_2<\ldots<x_n\}}\\

&=&
\sum_{\pi\in\mathfrak{S}: \pi(\varrho(i))=i, \forall i\in I_k)}
\frac{(\ind_{\{v\leq \mathfrak{p}(\varrho)_{j+1}\}}\zeta_\varrho u^M_v)_{\varrho(I_j)}(\mathbf{x}_{I_j})}{(\ind_{\{v\leq \mathfrak{p}(\varrho)_{j+1}\}}\zeta_\varrho a_{v-})_{\varrho(I_j)}(\mathbf{x}_{I_j})}a_{v-}(\pi(\mathbf{x}))\ind_{\{x_1<x_2<\ldots<x_n\}}\\

&=&
\frac{(\ind_{\{v\leq \mathfrak{p}(\varrho)_{j+1}\}}\zeta_\varrho u^M_v)_{\varrho(I_j)}(\mathbf{x}_{I_j})}{(\ind_{\{v\leq \mathfrak{p}(\varrho)_{j+1}\}}\zeta_\varrho a_{v-})_{\varrho(I_j)}(\mathbf{x}_{I_j})}\ind_{\{x_1<x_2<\ldots<x_n\}}
\sum_{\pi\in\mathfrak{S}: \pi(\varrho(i))=i \forall i\in I_k)}a_{v-}(\pi(\mathbf{x}))\\

&=&
\frac{(\ind_{\{v\leq \mathfrak{p}(\varrho)_{j+1}\}}\zeta_\varrho u^M_v)_{\varrho(I_j)}(\mathbf{x}_{I_j})}{(\ind_{\{v\leq \mathfrak{p}(\varrho)_{j+1}\}}\zeta_\varrho a_{v-})_{\varrho(I_j)}(\mathbf{x}_{I_j})}
\widetilde{a}^\varrho_{v-}.
\dce
$$
Applying Lemma \ref{attr2} the dual predictable projection $\wt N^\varrho_-\ast\psi(\wh V^\varrho)_t$ is given by$$
\sum_{j=0}^k\int_0^{t}\ind_{\{\tau_{(j)}< v\leq\tau_{(j+1)}\nmid\tau_{(k)}\}}
\frac{(\ind_{\{v\leq \mathfrak{p}(\varrho)_{j+1}\}}\zeta_\varrho u^M_v)_{\varrho(I_j)}(\overline{\btau}_{j}))}{(\ind_{\{v\leq \mathfrak{p}(\varrho)_{j+1}\}}\zeta_\varrho a_{v-})_{\varrho(I_j)}(\overline{\btau}_{j}))}
\frac{(\ind_{\{v\leq \mathfrak{p}_{j+1}\}}\widetilde{a}^\varrho_{v-})_{I_j}(\overline{\btau}_{j})}{(\ind_{\{v\leq \mathfrak{p}_{j+1}\}}\widetilde{a}_{v-})_{I_j}(\overline{\btau}_{j})}
d\cro{M,M}_v,
$$
for $t\in\rr_+$. 
\end{proof}

\brem
In the formula of the $\ff^{\overline{\btau}_{k}}$-drift of the bounded $\ff$-martingale $M$, we see clearly the term$$
\frac{
(\ind_{\{v\leq \mathfrak{p}(\varrho)_{j+1}\}}\zeta_\varrho u^M_v)_{\varrho(I_j)}(\overline{\btau}_{j}))}
{(\ind_{\{v\leq \mathfrak{p}(\varrho)_{j+1}\}}\zeta_\varrho a_{v-})_{\varrho(I_j)}(\overline{\btau}_{j}))}
$$
coming from the $\ff^{\btau_\varrho}$ decomposition of $M$, which appears in the formula with the weight$$
\frac{(\ind_{\{v\leq \mathfrak{p}_{j+1}\}}\widetilde{a}^\varrho_{v-})_{I_j}(\overline{\btau}_{j})}{(\ind_{\{v\leq \mathfrak{p}_{j+1}\}}\widetilde{a}_{v-})_{I_j}(\overline{\btau}_{j})}.
$$
\erem


\begin{thebibliography}{99} {\parskip = 3 pt



\bibitem{CJZ}
Callegaro, G. and Jeanblanc, M. and Zargari, B.: {\it Carthaginian
enlargement of filtrations} Forthcoming, ESAIM.



\bibitem{DM2}
Dellacherie, C. and Meyer, P.A.: {\it Probabilit\'es et potentiel}, Chapitres V-VIII: Th\'eorie des Martingales, Hermann (1992).

\bibitem{DMM}
Dellacherie, C., Maisonneuve, B and Meyer, P.A.: {\it Probabilit\'es et potentiel}, Chapitres XVII-XXIV: Processus de Markov (fin), Compl\'ements de
calcul stochastique, Hermann (1992).

\bibitem{EJJ} El-Karoui, N., Jeanblanc, M., and Jiao, Y. \emph{Modelling
several defaults}, Working paper. 2013.


\bibitem{GN}
Goutte, S and Ngoupeyou, A.: {\it Optimization problem and mean variance hedging on defaultable claims}. arXiv:1209.5953, Working paper, 2013.


\bibitem{HWY} He, S.W., Wang, J.G. and Yan. J.A.:
 {\it Semimartingale Theory and Stochastic Calculus.}
  Chinese edition, Science Press, Beijing, 1995.
%




\bibitem{Ja} Jacod ,~J.
{\it Calcul stochastique et probl\`emes de martingales}.
{\it S\'eminaire de Probabilit\'es XII, Lecture Notes in Mathematics 714}.
Springer,  Berlin Heidelberg New York, 1979.


\bibitem{JJ} Jacod, J.:
Grossissement initial, hypoth\`ese $(H')$ et th\'eor\`eme de Girsanov.
In: {\it Grossissements de Filtrations: Exemples et Applications. Lecture Notes in Mathematics 1118},
Th. Jeulin and M. Yor, eds. Springer,  Berlin Heidelberg New York, 1985, pp. 15--35.

\bibitem{J1} Jeulin,~T.:
Grossissement d'une filtration et applications.
In: {\it S\'eminaire de Probabilit\'es XIII, Lecture Notes in Mathematics 721},
Springer,  Berlin Heidelberg New York, 1979, pp. 574--609.

\bibitem{J2}
Jeulin,~T.: {\it Semi-martingales et grossissement d'une filtration}.
{\it Lecture Notes in Mathematics 833}. Springer,  Berlin Heidelberg New York, 1980.

\bibitem{JL}
Jeanblanc,~M. and Le Cam,~Y.: Progressive enlargement of filtration with initial times.
{\it Stochastic Processes and their Applications} 119 (2009), 2523--2543.



\bibitem{JKP}
Jiao, Y, Kharroubi, I and Pham, H
Optimal investment under multiple defaults risk: A BSDE-decomposition approach,
{\it Annals of Applied Probability}. Volume 23, 2 (2013), 455-491.

\bibitem{KP}
Kchia, Y. and Protter, P.: 
{\it On progressive filtration expansion with a process.} Working paper, 2011.


\bibitem{KLP}
Kchia,~Y., Larsson, M., and Protter,~P.:
Linking progressive and initial filtration expansions.
{\it Malliavin Calculus and Stochastic Analysis} A Festschrift in Honor of David Nualart Series: Springer Proceedings in Mathematics \& Statistics, Vol. 34  Viens, F., Feng, J., Hu, Y. and Nualart, E. eds. 2013

%




\bibitem{JY1} Jeulin,~T. and Yor,~M.:
Grossissement d'une filtration et semi-martingales: formules explicites.
In: {\it S\'eminaire de Probabilit\'es XII, Lecture Notes in Mathematics 649},
Springer,  Berlin Heidelberg New York, 1978, pp. 78--97.

\bibitem{M5} Meyer, P.A.:
Sur un th\'eor\`eme de J. Jacod
In: {\it S\'eminaire de Probabilit\'es VII, Lecture Notes in Mathematics 649.}
Springer-Verlag,  Berlin Heidelberg New York, 1978, 57--60.







\bibitem{RY}
Revuz, D. and Yor, M.:
{\it Continuous Martingales and Brownian Motion.} 3rd edition,
Springer, Berlin Heidelberg New York, 1999.

\bibitem{S} Song, S.: Optional splitting formula in progressively enlarged filtration,
to appear in ESAIM P\&S

\bibitem{St}Stricker, C.: Quasimartingales, martingales locales, et filtrations naturelles.
{\it Zeitschrift f\"ur Wahrscheinlichkeitstheorie and Verwandte Gebiete} 39 (1977), 55--63.

\bibitem{Y3} Yor,~M.:
Entropie d'une partition, et grossissement initial d'une filtration
In: {\it Grossissement de filtrations: exemples et applications. Lecture Notes 1118}.
Th. Jeulin and M. Yor, eds. Springer-Verlag, Berlin Heidelberg New York, 1980, 45--58.

}
\end{thebibliography}
\end{document}